\documentclass[11pt]{article}

\usepackage{amsmath,amssymb,amsfonts,amsbsy,mathrsfs}

\renewcommand{\subsection}{\subsubsection}

\newtheorem{theorem}{Theorem}[section]

\newtheorem{proposition}{Proposition}[section]

\newtheorem{remark}{Remark}[section]

\def\nt{|\hspace{-0.7pt}|\hspace{-0.7pt}|}

\begin{document}

\begin{center}
{\LARGE \bf On the well-posedness of a linearized plasma-vacuum}\\[6pt]
{\LARGE \bf interface problem in ideal compressible MHD}
\end{center}

\vspace*{7mm}

\centerline{\large Yuri Trakhinin}

\begin{center}
Sobolev Institute of Mathematics, Koptyug av. 4, 630090 Novosibirsk, Russia\\
e-mail: trakhin@math.nsc.ru
\end{center}

\vspace*{7mm}
\centerline{\bf \large Abstract}
\noindent We study the initial-boundary value problem resulting from the linearization of the plasma-vacuum interface problem in ideal compressible magnetohydrodynamics (MHD). We suppose that the plasma and the vacuum regions are unbounded domains and the plasma density does not go to zero continuously, but jumps. For the basic state upon which we perform linearization we find two cases of well-posedness of the ``frozen'' coefficient problem: the ``gas dynamical'' case and the ``purely MHD'' case. In the ``gas dynamical'' case we assume that the jump of the normal derivative of the total pressure is always negative. In the ``purely MHD''  case this condition can be violated but the plasma and the vacuum magnetic fields are assumed to be non-zero and non-parallel to each other everywhere on the interface. For this case we prove a basic a priori estimate in the anisotropic weighted Sobolev space $H^1_*$ for the variable coefficient problem.

\section{Introduction}
\label{s1}

Consider the equations of ideal compressible MHD with the gravitational field ${\cal G}\in \mathbb{R}^3$:
\begin{equation}
\left\{
\begin{array}{l}
\partial_t\rho  +{\rm div}\, (\rho {v} )=0,\\[3pt]
\partial_t(\rho {v} ) +{\rm div}\,(\rho{v}\otimes{v} -{H}\otimes{H} ) +
{\nabla}q=\rho{\cal G}, \\[3pt]
\partial_t{H} -{\nabla}\times ({v} {\times}{H})=0,\\[3pt]
\partial_t\bigl( \rho e +\frac{1}{2}|{H}|^2\bigr)+
{\rm div}\, \bigl((\rho e +p){v} +{H}{\times}({v}{\times}{H})\bigr)=\rho (v,{\cal G}),
\end{array}
\right.
\label{1}
\end{equation}
where $\rho$ denotes density, $v\in\mathbb{R}^3$ plasma velocity, $H \in\mathbb{R}^3$ magnetic field, $p=p(\rho,S )$ pressure, $q =p+\frac{1}{2}|{H} |^2$ total pressure, $S$ entropy, $e=E+\frac{1}{2}|{v}|^2$ total energy, and  $E=E(\rho,S )$ internal energy. With a state equation of gas, $p=p(\rho ,S)$, and the first principle of thermodynamics, \eqref{1} is a closed system. As the unknown we can fix, for example, the vector $ U =U (t, x )=(p, v,H, S)$.

System (\ref{1}) is supplemented by the divergent constraint
\begin{equation}
{\rm div}\, {H} =0
\label{2}
\end{equation}
on the initial data ${U} (0,{x} )={U}_0({x})$. As is known, taking into account \eqref{2}, we can easily symmetrize system \eqref{1} by rewriting it in the nonconservative form
\begin{equation}
\left\{
\begin{array}{l}
{\displaystyle\frac{1}{\rho c^2}}\,{\displaystyle\frac{{\rm d} p}{{\rm d}t} +{\rm div}\,{v} =0},\qquad
\rho\, {\displaystyle\frac{{\rm d}v}{{\rm d}t}-({H},\nabla ){H}+{\nabla}
q  =\rho {\cal G} },\\[9pt]
{\displaystyle\frac{{\rm d}{H}}{{\rm d}t} - ({H} ,\nabla ){v} +
{H}\,{\rm div}\,{v}=0},\qquad
{\displaystyle\frac{{\rm d} S}{{\rm d} t} =0},
\end{array}\right. \label{3}
\end{equation}
where $c^2=p_\rho(\rho ,S)$ is the square of the sound velocity and ${\rm d} /{\rm d} t =\partial_t+({v} ,{\nabla} )$ (by $(\ ,\ )$ we denote the scalar product). Equations (\ref{3}) read as the symmetric quasilinear system
\begin{equation}
\label{4}
A_0(U )\partial_tU+\sum_{j=1}^3A_j(U )\partial_jU+Q(U)=0,
\end{equation}
where $Q(U)=(0,-\rho{\cal G},0)$, $A_0= {\rm diag} \left(1/(\rho c^2),\rho ,\rho ,\rho ,
1,1,1,1\right)$,
\[
A_1=\left( \begin{array}{cccccccc} \frac{v_1}{\rho c^2}&1&0&0&0&0&0&0\\[6pt]
1&\rho v_1&0&0&0&{H_2}&{H_3}&0\\
0&0&\rho v_1&0&0&-{H_1}&0&0\\ 0&0&0&\rho v_1&0&0&-{H_1}&0\\
0&0&0&0&{v_1}&0&0&0\\
0&{H_2}&-{H_1}&0&0&{v_1}&0&0\\
0&{H_3}&0&-{H_1}&0&0&{v_1}&0\\ 0&0&0&0&0&0&0&v_1\\
\end{array} \right) ,
\]
\[
A_2=\left( \begin{array}{cccccccc} \frac{v_2}{\rho c^2}&0&1&0&0&0&0&0\\[6pt]
0&\rho v_2&0&0&-{H_2}&0&0&0\\ 1&0&\rho v_2&0&{H_1}&0&{H_3}&0\\ 0&0&0&\rho
v_2&0&0&-{H_2}&0\\
0&-{H_2}&{H_1}&0&{v_2}&0&0&0\\
0&0&0&0&0&{v_2}&0&0\\
0&0&{H_3}&-{H_2}&0&0&{v_2}&0\\0&0&0&0&0&0&0&v_2
\end{array} \right) ,
\]
\[
A_3=\left( \begin{array}{cccccccc} \frac{v_3}{\rho c^2}&0&0&1&0&0&0&0\\[6pt]
 0&\rho v_3&0&0&-{H_3}&0&0&0\\ 0&0&\rho v_3&0&0&-{H_3}&0&0\\ 1&0&0&\rho
v_3&{H_1}&{H_2}&0&0\\
0&-{H_3}&0&{H_1}&{v_3}&0&0&0\\
0&0&-{H_3}&{H_2}&0&{v_3}&0&0\\
0&0&0&0&0&0&{v_3}&0\\ 0&0&0&0&0&0&0&v_3
\end{array} \right)\; .
\]
System (\ref{4}) is symmetric hyperbolic if the the hyperbolicity condition $A_0>0$ holds:
\begin{equation}
\rho  >0,\quad p_\rho >0. \label{5}
\end{equation}

Plasma-vacuum interface problems for system \eqref{1} usually appear in the mathematical modeling of plasma confinement by magnetic fields. This subject was very popular in the 1950--70's, but most of theoretical studies were devoted to finding stability criteria of equilibrium states. The typical work in this direction is the classical paper of Bernstein et. al.
\cite{BFKK}. At the same time, according to our knowledge there are still no well-posedness results for full ({\it non-stationary}) plasma-vacuum models. Since \eqref{1} is a system of hyperbolic conservation laws which can produce shock waves and other types of strong discontinuities
(e.g., current-vortex sheets \cite{Tr}), it is natural to expect obtaining only local-in-time existence theorems.

The classical plasma-vacuum interface problem models confined plasmas in a closed vessel (see, e.g., \cite{Goed}). In this model the plasma is confined inside a perfectly conducting rigid wall and isolated from it by a vacuum region. Let $\Omega^+(t)$ and $\Omega ^-(t)$ are space-time domains occupied by the plasma and the vacuum respectively. That is, in the domain
$\Omega^+(t)$ we consider system \eqref{1} (or \eqref{4}) governing the motion of an ideal plasma and in the domain $\Omega^-(t)$  we have the elliptic (div-curl) system
\begin{equation}
\nabla \times \mathcal{H} =0,\qquad {\rm div}\, \mathcal{H}=0,\label{6}
\end{equation}
describing the vacuum magnetic field $\mathcal{H}\in\mathbb{R}^3$. Here, as in \cite{BFKK,Goed},
we consider so-called {\it pre-Maxwell dynamics}. That is, as usual in nonrelativistic MHD, we neglect the displacement current $(1/c)\,\partial_tE$, where $c$ is the speed of the light and $E$ is the electric field.

The boundary of the domain $\Omega^+(t)$ is a hypersurface $\Sigma (t)=\{F(t,x)=0\}$ that is the interface between plasma and vacuum. It is to be determined and moves with the velocity of plasma particles at the boundary:
\begin{equation}
\frac{{\rm d}F }{{\rm d} t}=0\quad \mbox{on}\ \Sigma (t)\label{7}
\end{equation} 
(for all $t\in [0,T]$). The plasma variable $U$ is connected with the vacuum magnetic field
$\mathcal{H}$  through the relations
\cite{BFKK,Goed}
\begin{equation}
[q]=0,\quad  (\mathcal{H}, N)=0,\quad (H,N)=0\quad \mbox{on}\ \Sigma (t),
\label{8}
\end{equation}
where $N=\nabla F$ and $[q]= q|_{\Sigma}-\frac{1}{2}|\mathcal{H}|^2_{|\Sigma}$. 
These relations together with \eqref{7} are the boundary conditions at the interface $\Sigma (t)$. At the perfectly conducting rigid wall $\Gamma$, that is the boundary of the vessel $\Omega=\Omega^-(t)\cup \Omega^+(t)$ and the exterior boundary of the vacuum region $\Omega^-(t)$, we have the boundary condition
\begin{equation}
(\mathcal{H}, n)=0\quad \mbox{on}\ \Gamma ,
\label{9}
\end{equation}
where $n$ is a normal vector to $\Gamma$.

From the mathematical point of view, a natural wish is to find conditions on the initial data
\begin{equation}
{U} (0,{x})={U}_0({x}),\quad {x}\in \Omega^{+} (0),\qquad
F(0,{x})=F_0({x}),\quad {x}\in\Sigma (0),
\label{10}
\end{equation}
\begin{equation}
\mathcal{H}(0,x)=
\mathcal{H}_0(x),\quad {x}\in \Omega^{-} (0),\label{11}
\end{equation}
providing the local-in-time existence and uniqueness of a solution $({U},\mathcal{H},F)$ of problem (\ref{1}), \eqref{6}--\eqref{11} in Sobolev spaces. Straightening the unknown interface
(see discussion below) and using the idea of the partition of unity, this complicated ``hyperbolic-elliptic'' free boundary problem could be splitted, roughly speaking, into two constituents:
problem \eqref{6}, \eqref{9}, \eqref{11} in the fixed domain $\Omega$ and problem
\eqref{1}, \eqref{6}--\eqref{8}, \eqref{10}, \eqref{11} with an unbounded domain $\Omega^-(t)$.
The first problem is reduced to the interior Neumann problem for the Laplace equation with a satisfied solvability condition by introducing the scalar potential $\Phi$, where $\nabla \Phi = \mathcal{H}(t,x)$. The second problem is our main interest in this paper and this problem is a natural generalization to MHD the free boundary problem for the compressible Euler equations with the ``vacuum'' boundary condition $p|_{\Sigma}=0$ (see \cite{Lind,Tr09}). For astrophysical plasmas this problem can be used for modeling the motion of a star when magnetic fields are taken into account.

As in \cite{Lind,Tr09}, we consider the case of liquid. This means that for problem \eqref{1}, \eqref{6}--\eqref{8}, \eqref{10}, \eqref{11} (with an unbounded domain $\Omega^-(t)$) the hyperbolicity conditions \eqref{5} are assumed to be satisfied in $\Omega^+$ up to the boundary $\Sigma$, i.e., the plasma density does not go to zero continuously, but jumps. At the same time, in the reality (e.g., for laboratory plasmas \cite{Goed}) the vacuum region is just a region of low enough density. That is, the assumption that the density is small but strictly positive at ${\Sigma}$ is quite reasonable.

Since the interface moves with the velocity of plasma particles at the boundary, at first sight
the passage to the Lagrangian coordinates to reduce the original problem to that in a fixed domain seems most natural. However, as, for example, for contact discontinuities in various models of fluid dynamics (e.g., for current-vortex sheets \cite{Tr}), this approach seems hardly realizable for problem \eqref{1}, \eqref{6}--\eqref{8}, \eqref{10}, \eqref{11}. Therefore, as in \cite{Tr09}, we will work in the Eulerian coordinates and for technical simplicity assume that the space-time domain $\Omega^+ (t)$ (the plasma region) is also unbounded and the interface $\Sigma (t)$ has the form of a graph:
$x_1=\varphi (t,x')$, $x'=(x_2,x_3)$. That is,  
\begin{equation}
\Omega^{\pm} (t)=\left\{x_1\gtrless \varphi (t,x')\right\}\label{12}
\end{equation}
and the function $\varphi (t,x')$ is to be determined. 

Now we can use a simple straightening of the unknown interface and reduce our problem to that in a half-space (see the next section). If,  however, the domain $\Omega^+ (t)$ is bounded and its initial boundary $\Sigma (0)$  is a compact co-dimension-1 surface in $\mathbb{R}^3$, as for shock waves, we can follow Majda's arguments \cite{Majda} (see also \cite[sect. 12.4.2]{BS}). More precisely, we can make ({\it locally in time}) a change of variables that sends all boundary locations $\Sigma (t)$ to the initial surface $\Sigma (0)$. We refer the reader to \cite{Majda,BS} for details of such a change of variables (see also \cite{Tr09} for further discussions).

Thus, we are finally interested in the following free boundary problem. We solve the symmetric {\it hyperbolic} system \eqref{4} (with assumption \eqref{5}) for $x_1>\varphi (t,x')$
and the {\it elliptic} system \eqref{6} for $x_1<\varphi (t,x')$. These systems are coupled through the boundary conditions \eqref{7}, \eqref{8} at the free boundary $x_1=\varphi (t,x')$. Moreover, we have the initial data \eqref{10}, \eqref{11} (with $F=x_1-\varphi (t,x')$) for $t=0$.

Actually, as for current-vortex sheets \cite{Tr}, we must regard the last boundary condition  in \eqref{8} as the restriction on the initial data \eqref{10}. More precisely, after straightening of the interface and in exactly the same manner as in \cite{Tr}, we can prove that a solution of \eqref{4}--\eqref{8}, \eqref{10}, \eqref{11} (if it exists for all $t\in [0,T]$) satisfies
\[
{\rm div}\, {H} =0 \quad \mbox{in}\ \Omega^+ (t)\quad \mbox{and}\quad (H,N)=0\quad \mbox{on}\ \Sigma (t)
\]
for all $t\in [0,T]$, if the latter was satisfied at $t=0$, i.e., for the initial data \eqref{10}. In particular, the fulfillment of ${\rm div}\, {H} =0$ implies that systems \eqref{1} and \eqref{4} are equivalent on solutions of problem \eqref{4}--\eqref{8}, \eqref{10}, \eqref{11}.

In the next section we first reduce problem \eqref{4}--\eqref{8}, \eqref{10}, \eqref{11} to that in the half-space $\mathbb{R}^3_+=\{x_1>0,\ x'\in \mathbb{R}^2\}$ and then linearize it about a basic state (``unperturbed flow''). For the basic state we consider two cases 
for which we can prove a priori estimates for the linearized problem with ``frozen'' (constant)
coefficients. In the first case, for the basic state we require the fulfillment of the condition
\begin{equation}
\left[\frac{\partial q}{\partial \mathcal{N}}\right]= -[\partial_1q]\leq -\epsilon <0,
\label{13}
\end{equation}
where $\mathcal{N}= (-1,0,0)$ is the outward normal to the boundary of $\mathbb{R}^3_+$, and
$[\partial_1q]=(\partial_1q)|_{x_1=0}-(\mathcal{H},\partial_1\mathcal{H})|_{x_1=0}.$
Since condition \eqref{13} is the counterpart of the natural physical condition  $\partial p/\partial N\leq -\epsilon <0$ in gas dynamics (see \cite{Lind,Tr09} and references therein), we call this case {\it ``gas dynamical.''} In the second case, we assume that the plasma and the vacuum tangential magnetic fields $(0,H_2,H_3)$ and $(0,\mathcal{H}_2,\mathcal{H}_3)$ are non-zero and non-parallel to each other everywhere on the straightened interface:
\begin{equation}
\left|(H_2\mathcal{H}_3-H_3\mathcal{H}_2)|_{x_1=0}\right|\geq \epsilon_1 >0.
\label{14}
\end{equation}
We call this case {\it ``purely MHD''} because for it the physical condition \eqref{13} can be violated, i.e., the magnetic field plays a stabilizing role.

From the mathematical point of view, the principal difference between the above cases is that for the ``purely MHD'' case the symbol associated to the interface is elliptic (see Section 4) and for the ``gas dynamical'' case this symbol can be non-elliptic. We suppose that one can prove a local-in-time existence and uniqueness theorem in Sobolev spaces for the original nonlinear problem for both of these cases, i.e., ``good'' initial data for the nonlinear problem reduced to that in the half-space $\mathbb{R}^3_+$ should satisfy either \eqref{13} or \eqref{14}. In this paper, we manage however to prove an a priori estimate for the variable coefficient linearized  problem only for the ``purely MHD'' case (see Section 5 for further discussions of open problems). 

The a priori estimate for the variable coefficient linearized problem that we derive for the ``purely MHD'' case can be considered as a first necessary step in proving the local-in-time existence for the original nonlinear interface problem by a suitable Nash-Moser-type iteration scheme. We plan to use the Nash-Moser method (as in \cite{Lind,Tr09}) because in this a priori estimate we have a loss of derivatives. Moreover, the additional difficulty is connected with the fact that the interface is a characteristic boundary for the hyperbolic system \eqref{4}.
This implies a natural loss of control on derivatives in the normal direction that cannot be compensated in MHD (unlike the situation in gas dynamics \cite{Sch,Tr09}). Therefore, the natural functional setting is provided by the anisotropic weighted Sobolev spaces $H^m_*$ (see \cite{YM,OSY,Sec,MST} and the next section for their definition). In this paper we prove our basic a priori estimate in the space $H^1_*$.

The rest of the paper is organized as follows. In Section 2, we obtain the linearized problem and formulate main results. In Section 3, for the constant coefficient linearized problem we derive a basic a priori $L_2$ estimate for the ``gas dynamical'' case. In Section 4, for the ``purely MHD'' case we prove an a priori estimate in $H^1_*$ for the variable coefficient problem. Section 5 is devoted to the discussion of open problems and contains concluding remarks. 

\section{Linearized problem and main results}
\label{s2}

\subsection{Reduction to a fixed domain}
\label{s2.1}

Let us first rewrite  the boundary conditions \eqref{7}, \eqref{8} for the unbounded domains
\eqref{12}:
\begin{equation}
\partial_t\varphi =v_N,\quad [q]=0,\quad \mathcal{H}_N=0,\quad H_N=0 \quad \mbox{on}\ \Sigma (t),\label{15}
\end{equation}
where $v_N=(v,N)$, $\mathcal{H}_N=(\mathcal{H},N)$, $H_N=(H,N)$, $N=(1,-\partial_2\varphi ,-\partial_2\varphi )$.

We straighten the interface $\Sigma$ by using the same change of independent variable as in \cite{Tr} (see also \cite{Tr09}). That is, the unknowns $U$ and $\mathcal{H}$ being smooth in $\Omega^{\pm}(t)$ are replaced by the vector-functions
\[
\widetilde{U}(t,x ):= {U}(t,\Phi^+ (t,x),x'),\quad 
\widetilde{\mathcal{H}}(t,x ):= \mathcal{H}(t,\Phi^- (t,x),x'),
\]
which are smooth in the half-space $\mathbb{R}^3_+$, where
\[
\Phi^{\pm}(t,x ):= \pm x_1+\Psi^{\pm}(t,x ),\quad \Psi^{\pm}(t,x ):= \chi (\pm x_1)\varphi (t,x').
\]
and $\chi\in C^{\infty}_0(\mathbb{R})$ equals to 1 on $[-1,1]$, and $\|\chi'\|_{L_{\infty}(\mathbb{R})}<1/2$. Here, as in \cite{Tr,Tr09}, we use the cut-off function $\chi$ to avoid assumptions about compact support of the initial data in our (future) nonlinear existence theorem. The above change of variable is admissible if $\partial_1\Phi^{\pm}\neq 0$. The latter is guaranteed, namely, the inequalities $\partial_1\Phi^+> 0$ and $\partial_1\Phi^-< 0$ are fulfilled, if we consider solutions for which $\|\varphi\|_{L_{\infty}([0,T]\times\mathbb{R}^2)}\leq 1$. This holds if,
without loss of generality, we consider the initial data satisfying $\|\varphi_0\|_{L_{\infty}(\mathbb{R}^2)}\leq 1/2$, and the time $T$ in our existence theorem is sufficiently small.

Dropping for convenience tildes in $\widetilde{U}$ and $\widetilde{\mathcal{H}}$, we reduce \eqref{4}, \eqref{6}, \eqref{15}, \eqref{10}, \eqref{11}  to the initial boundary value problem
\begin{equation}
\mathbb{P}(U,\Psi^+)=0,\quad \mathbb{V}(\mathcal{H},\Psi^-)=0\quad\mbox{in}\ [0,T]\times \mathbb{R}^3_+,\label{16}
\end{equation}
\begin{equation}
\mathbb{B}(U,\mathcal{H},\varphi )=0\quad\mbox{on}\ [0,T]\times\{x_1=0\}\times\mathbb{R}^{2},\label{17}
\end{equation}
\begin{equation}
(U,\mathcal{H})|_{t=0}=(U_0,\mathcal{H}_0)\quad\mbox{in}\ \mathbb{R}^3_+,\qquad \varphi|_{t=0}=\varphi_0\quad \mbox{in}\ \mathbb{R}^{2},\label{18}
\end{equation}
where $\mathbb{P}(U,\Psi^+)=P(U,\Psi^+)U$, 
\[
P(U,\Psi^+)=A_0(U)\partial_t +\widetilde{A}_1(U,\Psi^+)\partial_1+A_2(U )\partial_2+A_3(U )\partial_3,
\]
\[
\widetilde{A}_1(U,\Psi^+ )=\frac{1}{\partial_1\Phi^+}\Bigl(
A_1(U )-A_0(U)\partial_t\Psi^+-\sum_{k=2}^3A_k(U)\partial_k\Psi^+ \Bigr),
\]
\[
\partial_1\Phi^{\pm}=\pm 1 +\partial_1\Psi^{\pm},\qquad 
\mathbb{V}(\mathcal{H},\Psi^-)=\left(
\begin{array}{c}
\nabla\times \mathfrak{H}\\
{\rm div}\,\mathfrak{h}
\end{array}
\right),
\]
\[
\mathfrak{H}=(\mathcal{H}_1\partial_1\Phi^-,\mathcal{H}_{\tau_2},\mathcal{H}_{\tau_3}),\quad
\mathfrak{h}=(\mathcal{H}_{n},\mathcal{H}_2\partial_1\Phi^{-},\mathcal{H}_3\partial_1\Phi^{-}),
\]
\[
\mathcal{H}_{n}=\mathcal{H}_1-\mathcal{H}_2\partial_2\Psi^{-}-\mathcal{H}_3\partial_3\Psi^{-},\quad
\mathcal{H}_{\tau_i}=\mathcal{H}_1\partial_i\Psi^{-}+\mathcal{H}_i,\quad i=2,3,
\]
\[
\mathbb{B}(U,\mathcal{H},\varphi )=\left(
\begin{array}{c}
\partial_t\varphi -v_N \\ {[}q{]} \\ \mathcal{H}_N
\end{array}
\right),\quad v_{N}=v_1- v_2\partial_2\varphi - v_3\partial_3\varphi ,
\]
\[
[q]=q_{|x_1=0}-\frac{1}{2}|\mathcal{H}|^2_{x_1=0},\quad \mathcal{H}_N=
\mathcal{H}_1- \mathcal{H}_2\partial_2\varphi - \mathcal{H}_3\partial_3\varphi .
\]

In the MHD equations $\mathbb{P}(U,\Psi^+)=0$ we dropped, without loss of generality, the lower-order term $Q (U)$ responsible for the gravity. This term is only important for the case of {\it unbounded} domains for a correct configuration in Sobolev spaces of initial data satisfying conditions like \eqref{13} (see \cite{Tr09}). We just drop this term because it plays no role in our forthcoming linear analysis. We also did not include in our problem the equation
\begin{equation}
{\rm div}\, h=0\quad\mbox{in}\ [0,T]\times \mathbb{R}^3_+,\label{19}
\end{equation}
where $h=(H_{n},H_2\partial_1\Phi^+,H_3\partial_1\Phi^+)$,
$H_{n}=H_1-H_2\partial_2\Psi^+-H_3\partial_3\Psi^+$,
and the boundary condition
\begin{equation}
H_{N}=0\quad\mbox{on}\ [0,T]\times\{x_1=0\}\times\mathbb{R}^{2}\label{20}
\end{equation} 
because they are just restrictions on the initial data \eqref{18} (see Section 1). More precisely, referring to \cite{Tr} for the proof, we have the following proposition.

\begin{proposition}
Let the initial data \eqref{18} satisfy \eqref{19} and \eqref{20}.
If problem \eqref{16}--\eqref{18} has a solution $(U,\mathcal{H},\varphi )$, then this solution satisfies \eqref{19} and \eqref{20} for all $t\in [0,T]$. 
\label{p1}
\end{proposition}

Note that Proposition \ref{p1} stays valid if in \eqref{16} we replace system $\mathbb{P}(U,\Psi^+)=0$
by system \eqref{1} in the straightened variables. This means that these systems are equivalent on solutions of our plasma-vacuum interface problem and we may justifiably replace 
the conservation laws \eqref{1} by their nonconservative form \eqref{4}. 

\subsection{Basic state}
\label{s2.2}
Let
\begin{equation}
(\widehat{U}(t,x ),\widehat{\mathcal{H}}(t,x ),\hat{\varphi}(t,{x}')
\label{21}
\end{equation}
be a given sufficiently smooth vector-function 
with $\widehat{U}=(\hat{p},\hat{v},\widehat{H},\widehat{S})$ and
\begin{equation}
\|\widehat{U}\|_{W^2_{\infty}(\Omega_T)}+
\|\widehat{\mathcal{H}}\|_{W^2_{\infty}(\Omega_T)}+
\|\partial_1\widehat{U}\|_{W^2_{\infty}(\Omega_T)}+
\|\hat{\varphi}\|_{W^3_{\infty}(\partial\Omega_T)} \leq K,
\label{22}
\end{equation}
where $K>0$ is a constant and
\[
\Omega_T:= (-\infty, T]\times\mathbb{R}^3_+,\quad \partial\Omega_T:=(-\infty ,T]\times\{x_1=0\}\times\mathbb{R}^{2}.
\] 
If the basic state \eqref{21} upon which we shall linearize problem \eqref{16}--\eqref{18} is a solution of this problem (its existence should be proved), then it is natural to call it unperturbed flow. The trivial example of the unperturbed flow is the constant solution $(\overline{U},\overline{\mathcal{H}},0)$, where $\overline{U}\in \mathbb{R}^{8}$  and
$\overline{\mathcal{H}}\in \mathbb{R}^{3}$ are constant vectors.

We assume that the basic state \eqref{21} satisfies the hyperbolicity condition \eqref{5} in $\overline{\Omega_T}$,
\begin{equation}
\rho (\hat{p},\widehat{S}) >0,\quad \rho_p(\hat{p},\widehat{S}) >0 ,
\label{23}
\end{equation}
the first and the third boundary conditions in \eqref{17} on $\partial\Omega_T$,
\begin{equation}
\partial_t\hat{\varphi}-\hat{v}_{N}|_{x_1=0}=0,\quad \widehat{\mathcal{H}}_N|_{x_1=0}=0,
\label{24}
\end{equation}
and system $\mathbb{V}(\mathcal{H},\Psi^-)=0$ in ${\Omega_T}$,
\begin{equation}
\nabla\times \widehat{\mathfrak{H}}=0,\quad {\rm div}\,\hat{\mathfrak{h}}=0,
\label{25}
\end{equation}
where the ``hat'' values are determined like corresponding values for $(U,\mathcal{H},\varphi)$, e.g.,
\[
\widehat{\Phi}^{\pm}(t,x )=\pm x_1 +\widehat{\Psi}^{\pm}(t,x ),\quad
\widehat{\Psi}^{\pm}(t,x )=\chi(\pm x_1)\hat{\varphi}(t,x'),
\]
\[
\hat{v}_{N}=\hat{v}_1- \hat{v}_2\partial_2\hat{\varphi}- \hat{v}_3\partial_3\hat{\varphi},\quad
\widehat{\mathfrak{H}}=(\widehat{\mathcal{H}}_1\partial_1\widehat{\Phi}^-,
\widehat{\mathcal{H}}_{\tau_2},\widehat{\mathcal{H}}_{\tau_3}).
\]
Moreover, without loss of generality we assume that $\|\hat{\varphi}\|_{L_{\infty}(\partial\Omega_T)}<1$. This implies
\[
\partial_1\widehat{\Phi}^+\geq 1/2,\quad \partial_1\widehat{\Phi}^-\leq -  1/2.
\]
Note that (\ref{22}) yields
\[
\|\widehat{W}\|_{W^2_{\infty}(\Omega_T)} \leq C(K),
\]
where $\widehat{W}:=(\widehat{U},\partial_1\widehat{U},\widehat{\mathcal{H}},\nabla_{t,x}\widehat{\Psi}^+,
\nabla_{t,x}\widehat{\Psi}^-)$, $\nabla_{t,x}=(\partial_t, \nabla )$, and $C=C(K)>0$ is a constant depending on $K$.

\begin{remark}{\rm
Assumptions \eqref{23}--\eqref{25} are nonlinear constraints on the basic state. We will really need them while deriving a priori estimates for the linearized problem. In the forthcoming nonlinear analysis we plan to use the Nash-Moser method. As in \cite{Tr,Tr09}, the Nash-Moser procedure will be not completely standard. Namely, at each $n$th Nash-Moser iteration step we will have to construct an intermediate state $(U_{n+1/2},\mathcal{H}_{n+1/2},\varphi_{n+1/2})$ satisfying constraints  \eqref{23}--\eqref{25}. Without assumption \eqref{25} such an intermediate state can be constructed in exactly the same manner as in \cite{Tr,Tr09}. Assumption \eqref{25} does not however cause additional difficulties because for given $\hat{\varphi}$ it forms together with the last condition in \eqref{24} a boundary value problem reduced to the Neumann problem for the Laplace equation in the half-space. We omit corresponding arguments and postpone them to the nonlinear analysis.}
\label{r2.1}
\end{remark}

Later on, for the linearized problem we will need equations associated to the nonlinear constraints (\ref{19}) and (\ref{20}). However, to deduce them it is not enough that these constraints are satisfied by the basic state (\ref{21}). As in \cite{Tr}, we need actually that the equation for $H$ itself contained in system $\mathbb{P}(U,\Psi^+)=0$ is fulfilled for (\ref{21}):
\begin{equation}
\partial_t\widehat{H}+\frac{1}{\partial_1\widehat{\Phi}^+}\left\{ (\hat{w} ,\nabla )
\widehat{H} - (\hat{h} ,\nabla ) \hat{v} + \widehat{H}{\rm div}\,\hat{u}\right\} =0,
\label{26}
\end{equation}
where 
\[
\hat{u}=(\hat{v}_{n},\hat{v}_2\partial_1\widehat{\Phi}^+,\hat{v}_3\partial_1\widehat{\Phi}^+),\quad \hat{v}_{n}=\hat{v}_1-\hat{v}_2\partial_2\widehat{\Psi}^{+}-\hat{v}_3\partial_3\widehat{\Psi}^{+},
\]
and $\hat{w}=\hat{u}-(\partial_t\widehat{\Psi}^{+},0,0)$. Assume that (\ref{21}) satisfies (\ref{26}). Then, it follows from the proof of Proposition \ref{p1} (see \cite{Tr}) that constraints (\ref{19}) and (\ref{20}) are satisfied for the basic state (\ref{21}) if they are true for it at $t=0$. That is, without loss of generality we may suppose that (\ref{21}) satisfies
(\ref{19}) and (\ref{20}):
\begin{equation}
{\rm div}\,\hat{h}=0,\quad \widehat{H}_{N}|_{x_1=0}=0.
\label{27}
\end{equation}
Thus, for the basic state we require the fulfillment of conditions \eqref{22}--\eqref{27}.

\subsection{Linearized problem}
\label{s2.3}

The linearized equations for (\ref{16}), (\ref{17}) read:
\[
\mathbb{P}'(\widehat{U},\widehat{\Psi}^{+})(\delta U,\delta\Psi^{+}):=
\frac{\rm d}{{\rm d}\varepsilon}\mathbb{P}(U_{\varepsilon},\Psi_{\varepsilon}^{+})|_{\varepsilon =0}=f
\quad \mbox{in}\ \Omega_T,
\]
\[
\mathbb{V}'(\widehat{\mathcal{H}},\widehat{\Psi}^{-})(\delta \mathcal{H},\delta\Psi^{-}):=
\frac{\rm d}{{\rm d}\varepsilon}\mathbb{V}(\mathcal{H}_{\varepsilon},\Psi_{\varepsilon}^{-})|_{\varepsilon =0}=\mathcal{F}
\quad \mbox{in}\ \Omega_T,
\]
\[
\mathbb{B}'(\widehat{U},\widehat{\mathcal{H}},\hat{\varphi})(\delta U,\delta \mathcal{H},\delta \varphi ):=
\frac{\rm d}{{\rm d}\varepsilon}\mathbb{B}(U_{\varepsilon},\mathcal{H}_{\varepsilon},\varphi_{\varepsilon})|_{\varepsilon =0}={g}
\quad \mbox{on}\ \partial\Omega_T,
\]
where $U_{\varepsilon}=\widehat{U}+ \varepsilon\,\delta U$, $\mathcal{H}_{\varepsilon}=
\widehat{\mathcal{H}}+\varepsilon\,\delta \mathcal{H}$,
$\varphi_{\varepsilon}=\hat{\varphi}+ \varepsilon\,\delta \varphi$, and
\[
\Psi_{\varepsilon}^{\pm}(t,x ):=\chi (\pm x_1)\varphi_{\varepsilon}(t,x'),\quad
\Phi_{\varepsilon}^{\pm}(t,x ):=\pm x_1+\Psi_{\varepsilon}^{\pm}(t,x ),
\]
\[
\delta\Psi^{\pm}(t,x ):=\chi (\pm x_1)\delta \varphi (t,x ).
\]
Here we introduce the source terms $f=(f_1,\ldots ,f_8)$, $\mathcal{F}=(\mathcal{F}_1,\ldots ,
\mathcal{F}_4)$, and $g=(g_1,g_2,g_3)$ to make the interior equations and the boundary conditions inhomogeneous. 

We compute the exact form of the linearized equations (below we drop $\delta$): 
\[
\mathbb{P}'(\widehat{U},\widehat{\Psi}^+)(U,\Psi^+)
=
P(\widehat{U},\widehat{\Psi}^+)U +{\cal C}(\widehat{U},\widehat{\Psi}^+)
U -   \bigl\{L(\widehat{U},\widehat{\Psi}^+)\Psi^+\bigr\}\frac{\partial_1\widehat{U}}{\partial_1\widehat{\Phi}^+}
=f,
\]
\[
\mathbb{V}'(\widehat{\mathcal{H}},\widehat{\Psi}^{-})(\mathcal{H},\Psi^{-})=
\mathbb{V}(\mathcal{H},\widehat{\Psi}^{-})+
\left(\begin{array}{c}
\nabla\widehat{\mathcal{H}}_1\times\nabla\Psi^-\\[3pt]
\Bigl(\nabla \times \left(\begin{array}{c} 0 \\ -\widehat{\mathcal{H}}_3 \\
\widehat{\mathcal{H}}_2 \end{array} \right) ,\nabla\Psi^-\Bigr)
\end{array}
\right)=\mathcal{F},
\]
\[
\mathbb{B}'(\widehat{U},\widehat{\mathcal{H}},\hat{\varphi})(U,\mathcal{H},\varphi )=
\left(
\begin{array}{c}
\partial_t\varphi +\hat{v}_2\partial_2\varphi+\hat{v}_3\partial_3\varphi -v_{N}\\[3pt]
q-(\widehat{\mathcal{H}},\mathcal{H})\\[3pt]
\mathcal{H}_N-\widehat{\mathcal{H}}_2\partial_2\varphi -\widehat{\mathcal{H}}_3\partial_3\varphi
\end{array}
\right)_{|x_1=0}=g,
\]
where $q:=p+ (\widehat{H},H)$, $v_{N}:= v_1-v_2\partial_2\hat{\varphi}-v_3\partial_3\hat{\varphi}$, and the matrix
${\cal C}(\widehat{U},\widehat{\Psi}^+)$ is determined as follows:
\[
\begin{array}{r}
{\cal C}(\widehat{U},\widehat{\Psi}^{+})Y
= (Y ,\nabla_yA_0(\widehat{U} ))\partial_t\widehat{U}
 +(Y ,\nabla_y\widetilde{A}_1(\widehat{U},\widehat{\Psi}^{+}))\partial_1\widehat{U}
 \\[6pt]
+ (Y ,\nabla_yA_2(\widehat{U} ))\partial_2\widehat{U}
+ (Y ,\nabla_yA_3(\widehat{U} ))\partial_3\widehat{U},
\end{array}
\]
\[
(Y ,\nabla_y A(\widehat{U})):=\sum_{i=1}^8y_i\left.\left(\frac{\partial A (Y )}{
\partial y_i}\right|_{Y =\widehat{U}}\right),\quad Y =(y_1,\ldots ,y_8).
\]

The differential operators $\mathbb{P}'(\widehat{U},\widehat{\Psi}^{+})$ and $\mathbb{V}'(\widehat{\mathcal{H}},\widehat{\Psi}^{-})$ are first-order operators in $\Psi^{+}$
and $\Psi^{-}$ respectively. As in \cite{Tr}, following Alinhac \cite{Al}, we introduce the ``good unknown''
\begin{equation}
\dot{U}:=U -\frac{\Psi^+}{\partial_1\widehat{\Phi}^+}\,\partial_1\widehat{U}
\label{28}
\end{equation}
for the hyperbolic system of linearized MHD equations. Similarly, we also introduce the ``good unknown''
\begin{equation}
\dot{\mathcal{H}}:=\mathcal{H} -\frac{\Psi^-}{\partial_1\widehat{\Phi}^-}\,\partial_1\widehat{\mathcal{H}}
\label{29}
\end{equation}
for the elliptic system for the perturbation of the vacuum magnetic field.
Taking into account assumptions \eqref{24} and \eqref{25} and omitting detailed calculations,
we rewrite our linearized equations in terms of the new unknowns \eqref{28} and \eqref{29}:
\begin{equation}
P(\widehat{U},\widehat{\Psi}^+)\dot{U} +{\cal C}(\widehat{U},\widehat{\Psi}^+)
\dot{U} - \frac{\Psi^+}{\partial_1\widehat{\Phi}^+}\,\partial_1\bigl\{\mathbb{L}
(\widehat{U},\widehat{\Psi}^+)\bigr\}=f,
\label{30}
\end{equation}
\begin{equation}
\mathbb{V}(\dot{\mathcal{H}},\widehat{\Psi}^{-})=\mathcal{F}.
\label{31}
\end{equation}
\begin{multline}
\mathbb{B}'(\widehat{U},\widehat{\mathcal{H}},\hat{\varphi})(\dot{U},\dot{\mathcal{H}},\varphi ):= \mathbb{B}'(\widehat{U},\widehat{\mathcal{H}},\hat{\varphi})(U,\mathcal{H},\varphi )\\[6pt] =
 \left(
\begin{array}{c}
\partial_t\varphi+\hat{v}_2\partial_2\varphi+\hat{v}_3\partial_3\varphi-\dot{v}_{N}-
\varphi\,\partial_1\hat{v}_{N}\\[3pt]
\dot{q}-(\widehat{\mathcal{H}},\dot{\mathcal{H}})+ [\partial_1\hat{q}]\varphi \\[3pt]
\dot{\mathcal{H}}_{N}-\partial_2\bigl(\widehat{\mathcal{H}}_2\varphi \bigr) -\partial_3\bigl(\widehat{\mathcal{H}}_3\varphi \bigr)
\end{array}\right)_{|x_1=0}=g,
\label{32}
\end{multline}
where $\dot{v}_{\rm N}=\dot{v}_1-\dot{v}_2\partial_2\hat{\varphi}-\dot{v}_3\partial_3\hat{\varphi}$,
$\dot{\mathcal{H}}_{N}=\dot{v}_1-\dot{\mathcal{H}}_2\partial_2\hat{\varphi}-\dot{\mathcal{H}}_3\partial_3\hat{\varphi}$, and
\[
[\partial_1\hat{q}]=(\partial_1\hat{q})|_{x_1=0}-(\widehat{\mathcal{H}},\partial_1\widehat{\mathcal{H}})|_{x_1=0}.
\]
We used the last equation in \eqref{25} taken at $x_1=0$ while writing down the last boundary condition in \eqref{32}.

As in \cite{Al,Tr,Tr09}, we drop the zeroth-order term in $\Psi^+$ in (\ref{30}) and consider the effective linear operators
\[
\mathbb{P}'_e(\widehat{U},\widehat{\Psi}^+)\dot{U} :=P(\widehat{U},\widehat{\Psi}^+)\dot{U} +{\cal C}(\widehat{U},\widehat{\Psi}^+)
\dot{U}=f.
\]
In the future nonlinear analysis the dropped term in (\ref{30}) should be considered as an error term at each Nash-Moser iteration step. 

Regarding system \eqref{31}, without loss of generality we may actually drop the source term $\mathcal{F}$. At first sight, we have to keep it because the nonlinear system $\mathbb{V}(\mathcal{H},\Psi^-)=0$  will produce errors in the Nash-Moser iteration scheme. That is, in the future nonlinear analysis we will have to go
outside the class of divergence-free irrotational fields. At the same time, it follows from the detailed analysis of an exact form of the accumulated errors for the elliptic system $\mathbb{V}(\mathcal{H},\Psi^-)=0$ and the boundary condition $\mathcal{H}_N|_{x_1=0}=0$ (corresponding arguments are omitted and postponed to the nonlinear analysis) that the source terms ${\cal F}$ and $g_3$ have the following special form:
\begin{equation}
\mathcal{F} =\left(
\begin{array}{c}
\nabla\times\mathfrak{B}\\ {\rm div}\,\mathfrak{b}
\end{array}
\right),\quad g_3=\mathfrak{b}_1|_{x_1=0},\label{33}
\end{equation}
where
\[
\mathfrak{B}=(b_1\partial_1\Phi^-,b_{\tau_2},b_{\tau_3}),\quad
\mathfrak{b}=(b_{n},b_2\partial_1\Phi^{-},b_3\partial_1\Phi^{-}),
\]
\[
b_{n}=b_1-b_2\partial_2\Psi^{-}-b_3\partial_3\Psi^{-},\quad
b_{\tau_i}=b_1\partial_i\Psi^{-}+b_i,\quad i=2,3,
\]
and $b(t,x) =(b_1,b_2,b_3)$ is a vector-function. Passing to the new unknown
\[
\dot{\mathcal{H}}'=\dot{\mathcal{H}}-b
\]
and omitting then the primes, in view of \eqref{33}, we get the homogeneous system
\[
\mathbb{V}(\dot{\mathcal{H}},\widehat{\Psi}^{-})=0
\]
and the last boundary condition in \eqref{32} becomes homogeneous ($g_3=0$).

We now write down the final form of our linearized problem for $(\dot{U},\dot{\mathcal{H}},\varphi )$:
\begin{equation}
\widehat{A}_0\partial_t\dot{U}+\sum_{j=1}^{3}\widehat{A}_j\partial_j\dot{U}+
\widehat{\cal C}\dot{U}=f \qquad \mbox{in}\ \Omega_T,\label{34}
\end{equation}
\begin{equation}
\left\{ \begin{array}{ll}
\partial_t\varphi=\dot{v}_{N}-\hat{v}_2\partial_2\varphi-\hat{v}_3\partial_3\varphi +
\varphi\,\partial_1\hat{v}_{N}+g_1, & \\[6pt]
\dot{q}=(\widehat{\mathcal{H}},\dot{\mathcal{H}})-  [ \partial_1\hat{q}] \varphi +g_2 & \mbox{on}\ \partial\Omega_T, \end{array}\right.\label{35}
\end{equation}
\begin{equation}
\nabla\times \dot{\mathfrak{H}}=0,\quad {\rm div}\,\dot{\mathfrak{h}}=0 \qquad \mbox{in}\ \Omega_T, \label{36}
\end{equation}
\begin{equation}
\dot{\mathcal{H}}_{N} =\partial_2\bigl(\widehat{\mathcal{H}}_2\varphi \bigr) +\partial_3\bigl(\widehat{\mathcal{H}}_3\varphi \bigr)\qquad \mbox{on}\ \partial\Omega_T,
\label{37}
\end{equation}
\begin{equation}
(\dot{U},\dot{\mathcal{H}},\varphi )=0\qquad \mbox{for}\ t<0,\label{38}
\end{equation}
where 
\[
\widehat{A}_{\alpha}=:{A}_{\alpha}(\widehat{U}),\quad \alpha =0,2,3,\quad
\widehat{A}_1=:\widetilde{A}_1(\widehat{U},\widehat{\Psi}^+),\quad 
\widehat{\cal C}:={\cal C}(\widehat{U},\widehat{\Psi}^+),
\]
\[
\dot{\mathfrak{H}}=(\dot{\mathcal{H}}_1\partial_1\widehat{\Phi}^-,\dot{\mathcal{H}}_{\tau_2},\dot{\mathcal{H}}_{\tau_3}),\quad
\dot{\mathfrak{h}}=(\dot{\mathcal{H}}_{n},\dot{\mathcal{H}}_2\partial_1\widehat{\Phi}^{-},\dot{\mathcal{H}}_3\partial_1\widehat{\Phi}^{-}),
\]
\[
\dot{\mathcal{H}}_{n}=\dot{\mathcal{H}}_1-\dot{\mathcal{H}}_2\partial_2\widehat{\Psi}^{-}-\dot{\mathcal{H}}_3\partial_3\widehat{\Psi}^{-},\quad
\dot{\mathcal{H}}_{\tau_i}=\dot{\mathcal{H}}_1\partial_i\widehat{\Psi}^{-}+\dot{\mathcal{H}}_i,\quad i=2,3.
\]
We assume that $f$ and $g=(g_1,g_2)$ vanish in the past and consider the case of zero initial data, which is the usual assumption. We postpone the case of nonzero initial data to the nonlinear analysis (construction of a so-called approximate solution).

\subsection{Basic a priori estimates}
\label{s2.4}

We first write down our basic a priori estimates for the case of constant (``frozen'') coefficients of problem \eqref{34}--\eqref{38}. Before formulating this result we give the definition of the anisotropic weighted Sobolev spaces $H^m_*$. Following \cite{YM,OSY,Sec,MST}, the functional space $H^m_*$ is defined as follows:
$$
H^{m}_*(\mathbb{R}^3_+ ):=\left\{ u\in L_2(\mathbb{R}^3_+) \ | \
\partial^{\alpha}_*\partial_1^k u\in L_2(\mathbb{R}^3_+ )\quad \mbox{if}\quad |\alpha
|+2k\leq m\, \right\},
$$
where $m\in\mathbb{N}$, $\partial^{\alpha}_*=(\sigma \partial_1)^{\alpha_1}\partial_2^{\alpha_2}
\partial_3^{\alpha_3}\,$, and  $\sigma (x_1)\in C^{\infty}(\mathbb{R}_+)$ is
a monotone increasing function such that $\sigma (x_1)=x_1$ in a neighborhood of
the origin and $\sigma (x_1)=1$ for $x_1$ large enough. The space $H^m_*(\mathbb{R}^3_+ )$ is normed by
\[
\|u\|_{m,*}^2= \sum_{|\alpha |+2k\leq m}
\|\partial^{\alpha}_*\partial_1^k u\|_{L_2(\mathbb{R}^3_+)}^2.
\]
We also define the space
\[
H^m_*(\Omega_T)=\bigcap_{k=0}^{m}H^k((-\infty ,T],H^{m-k}_*(\mathbb{R}^3_+ ))
\]
equipped with the norm
\[
[u]^2_{m,*,T}=\int_{-\infty}^T\nt u (t)\nt^2_{m,*}dt,
\quad \mbox{where}\quad
\nt u (t)\nt^2_{m,*}=\sum\limits_{j=0}^m\|\partial_t^ju(t)\|^2_{m-j,*}.
\]
Within this paper we use the space $H^m_*(\Omega_T)$ mainly for $m=1$. Clearly, 
the norm for $H^1_*(\Omega_T)$ reads
\[
[u]^2_{1,*,T}=\int_{\Omega_T}\left( u^2 +(\partial_tu)^2 +(\sigma\partial_1u)^2 
+(\partial_2u)^2+(\partial_3u)^2\right) dtdx.
\]

We are now in a position to state our main results.

\begin{theorem}
Let the basic state \eqref{21} satisfies assumptions \eqref{23}--\eqref{27}. Let the coefficients of problem \eqref{34}--\eqref{38} are ``frozen'', i.e., the coefficients of the interior equations \eqref{34}, \eqref{36} and the coefficients of the boundary conditions \eqref{35}, \eqref{37}  have been calculated at given points
$(t_*,x_1^*,x'_*)\in \Omega_T$ and $(t_*,x'_*)\in \partial\Omega_T$ respectively, in particular, the coefficient $[\partial_1\hat{q}]$ is a constant. Then, for the ``gas dynamical'' case \eqref{13},
\begin{equation}
[\partial_1\hat{q}] >0,\label{39}
\end{equation}
sufficiently smooth solutions $(\dot{U},\dot{\mathcal{H}},\varphi )$ of \eqref{34}--\eqref{38}
obey the estimate
\begin{equation}
\|\dot{U}\|_{L_2(\Omega_T)}+\|\dot{\mathcal{H}}\|_{L_2(\Omega_T)}+\|\varphi\|_{L_2(\partial\Omega_T)}\leq C\left\{ \|f\|_{L_2(\Omega_T)}+
\|g\|_{H^{1}(\partial\Omega_T)}\right\},
\label{40}
\end{equation}
where $C=C(T)>0$ is a constant independent of the data $(f,g)$.
\label{t2.2} 
\end{theorem}

For variable coefficients, we have not yet managed to derive an a priori estimate for the ``gas dynamical'' case, but for the ``purely MHD'' case we prove the following theorem.

\begin{theorem}
Let the basic state \eqref{21} satisfies assumptions \eqref{22}--\eqref{27}. 
Then, for the ``purely MHD'' case \eqref{14},
\begin{equation}
|(\widehat{H}_2\widehat{\mathcal{H}}_3-\widehat{H}_3\widehat{\mathcal{H}}_2)|_{x_1=0}|\geq \epsilon_1> 0,\label{41}
\end{equation}
sufficiently smooth solutions $(\dot{U},\dot{\mathcal{H}},\varphi )$ of problem \eqref{34}--\eqref{38} obey the estimate
\begin{equation}
[\dot{U}]_{1,*,T}+\|\dot{\mathcal{H}}\|_{H^{1}(\Omega_T)}+\|\varphi\|_{H^1(\partial\Omega_T)}\leq C\left\{ [f]_{2,*,T}+
\|g\|_{H^{2}(\partial\Omega_T)}\right\},
\label{42}
\end{equation}
where $C=C(K,T)>0$ is a constant independent of the data $(f,g)$.
\label{t2.3} 
\end{theorem}

\section{``Gas dynamical'' case}
\label{s3}

\subsection{Properties of problem \eqref{34}--\eqref{38}}

Before ``freezing'' coefficients we discuss some useful properties of the variable coefficient problem \eqref{34}--\eqref{38}. First of all, as for current-vortex sheets \cite{Tr}, we can prove the following proposition (see Appendix A in \cite{Tr} for the proof).

\begin{proposition}
Let the basic state \eqref{21} satisfies assumptions \eqref{22}--\eqref{27}. 
Then solutions of problem \eqref{34}--\eqref{38} satisfy
\begin{equation}
{\rm div}\,\dot{h}=r\quad\mbox{in}\ \Omega_T,
\label{43}
\end{equation}
\begin{equation}
\widehat{H}_2\partial_2\varphi +\widehat{H}_3\partial_3\varphi -\dot{H}_{N}-
\varphi\,\partial_1\widehat{H}_{N}=g_3\quad\mbox{on}\ \partial\Omega_T.
\label{44}
\end{equation}
Here
\[
\dot{h}=(\dot{H}_{n},\dot{H}_2\partial_1\widehat{\Phi}^+,\dot{H}_3\partial_1\widehat{\Phi}^+),\quad
\dot{H}_{n}=\dot{H}_1-\dot{H}_2\partial_2\widehat{\Psi}^+-\dot{H}_3
\partial_3\widehat{\Psi}^+
\]
($\dot{H}_{N}|_{x_1=0}=\dot{H}_{n}|_{x_1=0}$), the functions $r=
r(t,x )$ and $g_3= g_3(t,x')$, which vanish in the past, are determined by the source terms and the basic state as solutions to the linear inhomogeneous equations
\begin{equation}
\partial_t a+ \frac{1}{\partial_1\widehat{\Phi}^+}\left\{ (\hat{w} ,\nabla a) + a\,{\rm div}\,\hat{u}\right\}={\cal F}\quad\mbox{in}\ \Omega_T,
\label{45}
\end{equation}
\begin{equation}
\partial_t g_3 +\hat{v}_2\partial_2g_3+\hat{v}_3\partial_3g_3+
\left(\partial_2\hat{v}_2+\partial_3\hat{v}_3\right) g_3={\cal G}\quad \mbox{on}\ \partial\Omega_T,
\label{46}
\end{equation}
where $a=r/\partial_1\widehat{\Phi}^+,\quad {\cal F}=({\rm div}\,
{f}_{H})/\partial_1\widehat{\Phi}^+$,
\[
{f}_{H}=
(f_{n} ,f_6,f_7),\quad
f_{n}=f_5-f_6\partial_2\widehat{\Psi}^+-
f_7\partial_3\widehat{\Psi}^+,
\quad
{\cal G}=\bigl\{\partial_2\bigl(\widehat{H}_2g_1\bigr)+
\partial_3\bigl(\widehat{H}_3g_1\bigr)-f_{n}\bigr\}\bigr|_{x_1=0}.
\]
\label{p3.1}
\end{proposition}

It follows from the first condition in \eqref{24} that the interior equation \eqref{45} does not need a boundary condition because $\hat{w}_1|_{x_1=0}=0$. Therefore, from \eqref{45} we get
\begin{equation}
\|r\|_{L_2(\Omega_T)}\leq C\|f\|_{H^1(\Omega_T)}\leq C[f]_{2,*,T}.\label{47}
\end{equation}
Here an later on $C$ is a constant that can change from line to line, and sometimes we show the dependence of $C$ from another constants. In particular, in \eqref{47} the constant $C$ depends on $K$ and $T$. Using \eqref{46} and the trace theorem \cite{OSY} for the spaces $H^m_*$, we easily estimate:
\begin{equation}
\|g_3\|_{H^1(\partial\Omega_T)}\leq C\left\{\|g_1\|_{H^2(\partial\Omega_T)}+\|f_{|x_1=0}\|_{H^1(\partial\Omega_T)}\right\}
\leq C\{\|g\|_{H^2(\partial\Omega_T)}+ [f]_{2,*,T}\}.\label{48}
\end{equation}
%Estimates \eqref{47} and \eqref{48} cause additional loss of derivatives with respect to the %source terms in \eqref{42} (otherwise, we had the norms $[f]_{1,*,T}$ and %$\|g\|_{H^2(\partial\Omega_T)}$ in the right-hand side of \eqref{42}). This is a natural prize %for the non-inclusion of the divergent constraint and the boundary condition for the plasma %magnetic field in problem \eqref{34}--\eqref{38}

In view of the first condition in \eqref{24} and the second condition in \eqref{27}, the boundary matrix $\widehat{A}_1$ is singular at $x_1=0$, i.e., the plane $x_1=0$ is a characteristic boundary for the hyperbolic system \eqref{34} (exactly as for current-vortex sheets \cite{Tr}). Following \cite{Tr}, we introduce the new unknown 
\[
V=(\dot{q},\dot{v}_{n},\dot{v}_2,\dot{v}_3,\dot{H}_{n},\dot{H}_2,
\dot{H}_3,\dot{S})
\]
for separating ``characteristic'' and ``noncharacteristic'' unknowns.
We have $\dot{U}=JV$, with 
\[
J=\left(
\begin{array}{cccccccc}
1 & \; 0 & 0 & 0 & -\widehat{H}_1 & -\widehat{H}_{\tau_2} &
-\widehat{H}_{\tau_3} & 0 \\
0 & \; 1 & \;\partial_2\widehat{\Psi}^+ & \partial_3\widehat{\Psi}^+ & 0 & 0 & 0 & 0 \\
0 & \; 0 & 1 & 0 & 0 & 0 & 0 & 0 \\
0 &\; 0 & 0 & 1 & 0 & 0 & 0 & 0 \\
0 &\; 0 & 0 & 0 & 1 & \partial_2\widehat{\Psi}^+ & \partial_3\widehat{\Psi}^+ & 0 \\
0 &\; 0 & 0 & 0 & 0 & 1 & 0 & 0 \\
0 &\; 0 & 0 & 0 & 0 & 0 & 1 & 0 \\
0 &\; 0 & 0 & 0 & 0 & 0 & 0 & 1 \end{array}
\right),
\]
where $\widehat{H}_{\tau_i}$ ($=2,3$) are determined in the same way as $\widehat{\mathcal{H}}_{\tau_i}$ above. Then, system \eqref{34} is equivalently rewritten as
\begin{equation}
{\cal A}_0\partial_t{V}+ \sum_{k=1}^{3}{\cal A}_k\partial_k{V} +{\cal
A}_4V =F, \label{49}
\end{equation}
where
${\cal A}_{\alpha}=J^{\textsf{T}}\widehat{A}_{\alpha}J$ ($\alpha =\overline{0,3}$),
$F=J^{\textsf{T}}f$, and 
\[
{\cal A}_4=J^{\textsf{T}}\left\{ \widehat{\cal C}+\widehat{A}_0\partial_t J+\sum _{k=1}^3\widehat{A}_k\partial_k J\right\}.
\]

The boundary matrix ${\cal A}_1$ in system \eqref{49} has the form
\begin{equation}
{\cal A}_1={\cal A}+{\cal A}_{(0)},\quad
{\cal A}=\frac{1}{\partial_1\widehat{\Phi}^+}\,{\cal E}_{12},\quad {\cal A}_{(0)}|_{x_1=0}=0,
\label{50}
\end{equation}
where ${\cal E}_{ij}$ is the symmetric matrix which ($ij$)th and ($ji$)th elements equal to 1 and others are zero. The explicit form of ${\cal A}_{(0)}$ is of no interest, and it is only important that ${\cal A}_{(0)}|_{x_1=0}=0$. Therefore, 
the boundary matrix ${\cal A}_1$ on the boundary $x_1=0$ is of constant rank 2. That is, \eqref{49} is a symmetric hyperbolic system with characteristic boundary of constant multiplicity (in the sense of Rauch \cite{Rauch}). It is also noteworthy that because of
(\ref{43}) not only $\dot{q}$ and $\dot{v}_{n}$ but also $\dot{H}_{n}$  is a ``noncharacteristic'' unknown. For the ``noncharacteristic'' part of the vector $V$,
\begin{equation}
V_{n} =(\dot{q},\dot{v}_{n},\dot{H}_{n}),
\label{51}
\end{equation}
we expect to have a better control on the normal ($x_1$-) derivatives.

We now discuss the elliptic part of our problem. In \cite{BFKK}, the vacuum vector potential was used for the div-curl system. Unlike \cite{BFKK}, here we introduce the scalar potential $A$ for system \eqref{36}:
\begin{equation}
\dot{\mathfrak{H}}=\nabla A.
\label{52}
\end{equation}
Then, we get
\begin{equation}
\widetilde{\triangle} A =0,
\label{53}
\end{equation}
with
\[
\widetilde{\triangle} =\tilde{\partial}_1^2+\tilde{\partial}_2^2+\tilde{\partial}_3^2,\quad
\tilde{\partial}_1=\frac{1}{\partial_1\widehat{\Phi}^-}\,\partial_1,\quad
\tilde{\partial}_k=\partial_k-\frac{\partial_k\widehat{\Psi}^-}{\partial_1\widehat{\Phi}^-}\,\partial_1,\quad k=2,3.
\]
Passing to the ``original'' curvilinear coordinates $\tilde{x}_1=\widehat{\Phi}^-(t,x)$,
$\tilde{x}'=x'$, we could rewrite \eqref{53} and the boundary condition \eqref{37} as the Neumann problem for the Laplace equation for $A$ (with a satisfied solvability condition). 
However, we do not need do so. Moreover, we do not rewrite \eqref{37} in terms of the potential $A$. In fact we will only use relation \eqref{52} and it will be even more convenient for us to work with the equation 
\begin{equation}
{\rm div}\,\dot{\mathfrak{h}}=0
\label{54}
\end{equation} 
instead of equation \eqref{53}.

\subsection{Proof of Theorem \ref{t2.2}}

We first do not ``freeze'' coefficients and obtain an inequality for variable coefficients which will imply the a priori estimate \eqref{40} for ``frozen'' coefficients. 
By standard argument we obtain for the hyperbolic system \eqref{49}  the energy inequality
\begin{equation}
I(t)- 2\int_{\partial\Omega_t}\dot{q}\,\dot{v}_N|_{x_1=0}\,{\rm d}x'\,{\rm d}s\leq C 
\left( \| f\|^2_{L_2(\Omega_T)} +\int_0^tI(s)\,{\rm d}s\right),
\label{55}
\end{equation}
where $I(t)=\int_{\mathbb{R}^3_+}({\cal A}_0V,V)\,{\rm d}x$. It follows from the boundary
conditions \eqref{35} that
\begin{equation}
-2\dot{q}\,\dot{v}_N|_{x_1=0}=2\left([\partial_1\hat{q}]\varphi- g_2\right)\dot{v}_N|_{x_1=0}
-2(\widehat{\mathcal{H}},\dot{\mathcal{H}})\dot{v}_N|_{x_1=0},
\label{56}
\end{equation}
where
\begin{multline}
2\left([\partial_1\hat{q}]\varphi- g_2\right)\dot{v}_N|_{x_1=0} \\ =2([\partial_1\hat{q}]\varphi- g_2)(\partial_t\varphi+\hat{v}_2\partial_2\varphi+\hat{v}_3\partial_3\varphi -
\varphi\,\partial_1\hat{v}_{N}-g_1)|_{x_1=0}\\
=\partial_t\left\{ [\partial_1\hat{q}]\varphi^2
-2g_2\varphi\right\}+\partial_2\left\{ \hat{v}_2|_{x_1=0}\,[\partial_1\hat{q}]\varphi^2 -2\hat{v}_2|_{x_1=0}\,g_2\varphi \right\} \\ +
\partial_3\left\{ \hat{v}_3|_{x_1=0}\,[\partial_1\hat{q}]\varphi^2 -2\hat{v}_3|_{x_1=0}\,g_2\varphi  \right\}
+2g_1g_2
\\
-\bigl\{\partial_t[\partial_1\hat{q}]+\partial_2(\hat{v}_2[\partial_1\hat{q}])+\partial_3(\hat{v}_3[\partial_1\hat{q}])-2[\partial_1\hat{q}]\partial_1\hat{v}_{N}\bigr\}|_{x_1=0}\,\varphi^2
\\
+2\left\{ \partial_tg_2+\partial_2(\hat{v}_2g_2)+\partial_3(\hat{v}_3g_2)+g_2\partial_1\hat{v}_{N}-[\partial_1\hat{q}]g_1\right\}|_{x_1=0}\,\varphi .
\label{57}
\end{multline}

Assume that $[\partial_1\hat{q}] \geq \epsilon >0$ (this is the version of condition \eqref{39} for variable coefficients). Using the Young inequality, from \eqref{55}--\eqref{57} we obtain
\begin{multline}
I(t)+\frac{1}{2}\int_{\mathbb{R}^2}[\partial_1\hat{q}]\,\varphi^2\,{\rm d}x'
- 2\int_{\partial\Omega_t}(\widehat{\mathcal{H}},\dot{\mathcal{H}})\dot{v}_N|_{x_1=0}\,{\rm d}x'\,{\rm d}s
\\
\leq C(K,\epsilon ) \Bigl\{ \| f\|^2_{L_2(\Omega_T)}+\| g\|^2_{H^1(\partial\Omega_T)}
+\int_0^t\left(I(s)
+\|\varphi (s)\|^2_{L_2(\mathbb{R}^2)}\right)\,{\rm d}s\Bigr\}.\label{58}
\end{multline}
Regarding the boundary term $-2(\widehat{\mathcal{H}},\dot{\mathcal{H}})\dot{v}_N|_{x_1=0}$ in \eqref{58}, in view of \eqref{35}, \eqref{37}, \eqref{52}, and the second condition in \eqref{24}, we have
\begin{equation}
-2(\widehat{\mathcal{H}},\dot{\mathcal{H}})\dot{v}_N|_{x_1=0} =
-2(\widehat{\mathcal{H}}_2\partial_2A+\widehat{\mathcal{H}}_3\partial_3A)(\partial_t\varphi+\hat{v}_2\partial_2\varphi+\hat{v}_3\partial_3\varphi -
\varphi\,\partial_1\hat{v}_{N}-g_1)|_{x_1=0}.\label{59}
\end{equation}

Note that in the framework of the $L_2$ theory we are not able to treat the term $2g_1(\widehat{\mathcal{H}},\dot{\mathcal{H}})|_{x_1=0}$ contained in \eqref{59} directly. On the other hand, since in \eqref{58} we anyway lose one derivative from $g$, we can use the classical argument suggesting to reduce our problem to one with homogeneous boundary conditions by subtracting from the solution a more regular function. Namely, there exists $\widetilde{U}=(\tilde{p},\tilde{v},\widetilde{H},\tilde{S})\in H^{1}(\Omega_T)$ (or more precisely, $(\tilde{q},\tilde{v}_n)\in H^{1}(\Omega_T)$ and $(\tilde{v}_2,\tilde{v}_3,\widetilde{H},\tilde{S})\in H^{1}_*(\Omega_T)$) vanishing in the past such that 
\[
-\tilde{v}_N=g_1,\quad \tilde{q}=g_2\quad\mbox{on}\ \partial\Omega_T,
\]
where $\tilde{v}_n$, $\tilde{v}_{N}$, and $\tilde{q}$ are determined like corresponding values for $\dot{U}$.
If $\dot{U}=U^{\natural}+\widetilde{U}$, then $(U^{\natural},\dot{\mathcal{H}},\varphi )$ satisfies \eqref{34}--\eqref{38} with $g=0$ and
$f=f^{\natural}$, where $f^{\natural}= f-\mathbb{P}'_e(\widehat{U},\widehat{\Psi})\widetilde{U}$. That is, it is enough to prove estimate \eqref{40} with $g=0$. Without loss of generality, we will just assume that in \eqref{59} $g_1=0$.

Let us consider the term $-2\partial_t\varphi (\widehat{\mathcal{H}},\dot{\mathcal{H}})|_{x_1=0}$ contained in \eqref{59}. Integrating by parts and using \eqref{37} and \eqref{52}, we obtain
\begin{multline}
-2\int_{\partial\Omega_t}(\widehat{\mathcal{H}}_2\partial_2A+\widehat{\mathcal{H}}_3\partial_3A)
|_{x_1=0}\,\partial_t\varphi\,{\rm d}x'\,{\rm d}s \\ =-2\int_{\partial\Omega_t}\Bigl(
\nabla'A,\bigl(\partial_t(\varphi\widehat{\mathcal{H}}' )-\varphi\partial_t\widehat{\mathcal{H}}'\bigr)\Bigr)\Bigr|_{x_1=0}\,{\rm d}x'\,{\rm d}s \\
=2\int_{\partial\Omega_t}A\,\partial_t\dot{\mathcal{H}}_N|_{x_1=0}\,{\rm d}x'\,{\rm d}s
+2\int_{\partial\Omega_t}(\dot{\mathfrak{H}}',\partial_t\widehat{\mathcal{H}}')|_{x_1=0}\,\varphi\,{\rm d}x'\,{\rm d}s,
\label{60}
\end{multline}
where $\nabla'=(\partial_2,\partial_3)$,  $\widehat{\mathcal{H}}'=(\widehat{\mathcal{H}}_2,\widehat{\mathcal{H}}_3)$, and 
$\dot{\mathfrak{H}}'=(\dot{\mathcal{H}}_{\tau_2},\dot{\mathcal{H}}_{\tau_3})$.
It is noteworthy that for ``frozen'' coefficients the last integral in \eqref{60} disappears. Regarding the 
penultimate integral in \eqref{60}, taking into account \eqref{36} and \eqref{52}, we have
\begin{multline}
2\int_{\partial\Omega_t}A\,\partial_t\dot{\mathcal{H}}_N|_{x_1=0}\,{\rm d}x'\,{\rm d}s  =
J(t) -
2\int_{\partial\Omega_t}\dot{\mathcal{H}}_N\partial_tA|_{x_1=0}\,{\rm d}x'\,{\rm d}s \\ =
J(t)+2\int_{\Omega_t}\partial_1\bigl(\dot{\mathcal{H}}_n\partial_tA\bigr)\,{\rm d}x\,{\rm d}s 
 \\ =J(t)+
2\int_{\Omega_t}\Bigl\{ \dot{\mathcal{H}}_n\partial_t\bigl( \dot{\mathcal{H}}_1\partial_1\widehat{\Phi}^-\bigr) 
-\partial_tA \,{\rm div}_{x'}\bigl(\dot{\mathcal{H}}'\partial_1\widehat{\Phi}^-\bigr)\Bigr\}\,{\rm d}x\,{\rm d}s \\
= J(t) +2\int_{\Omega_t}\Bigl\{ \dot{\mathcal{H}}_n\partial_t\bigl( \dot{\mathcal{H}}_1\partial_1\widehat{\Phi}^-\bigr) \\ +
\partial_1\widehat{\Phi}^-\bigl( \dot{\mathcal{H}}', \partial_t(\dot{\mathcal{H}}_1\nabla'\widehat{\Psi}^- + \dot{\mathcal{H}}')\bigr)\Bigr\}\,{\rm d}x\,{\rm d}s 
= J(t) - K(t) +L(t),
\label{61}
\end{multline}
where $\dot{\mathcal{H}}'=(\dot{\mathcal{H}}_2,\dot{\mathcal{H}}_3)$, ${\rm div}_{x'}b':= \partial_2b_2+\partial_3b_3$ (with $b'=(b_2,b_3)\in \mathbb{R}^2$),
\[
J(t)=2\int_{\mathbb{R}^2}A\dot{\mathcal{H}}_N|_{x_1=0}\,{\rm d}x',\quad
K(t)= -\int_{\mathbb{R}^3_+}\partial_1\widehat{\Phi}^-|\dot{\mathcal{H}}|^2\,{\rm d}x,
\]
\[
L(t) =\int_{\Omega_t}\Bigl\{\partial_1\partial_t\widehat{\Psi}^-( \dot{\mathcal{H}}_1^2-|\dot{\mathcal{H}}'|^2)
+\bigl((\partial_1\widehat{\Phi}^-\nabla'\partial_t\widehat{\Psi}^--
\partial_1\partial_t\widehat{\Psi}^-\nabla'\widehat{\Psi}^-),\dot{\mathcal{H}}'\bigr)
\Bigr\}\,{\rm d}x\,{\rm d}s.
\]

Recall that $\partial_1\widehat{\Phi}^-\leq -1/2$, i.e., $K(t)\geq \frac{1}{2}\,\|\dot{\mathcal{H}}(t)\|^2_{L_2(\mathbb{R}^3_+)}$. 
For ``frozen'' coefficients $L(t)\equiv 0$, but even for the variable coefficients case we easily estimate $-L(t)$ from above by $C\|\dot{\mathcal{H}}\|^2_{L_2(\Omega_t)}$.
Multiplying equation \eqref{54} by the potential $A$, integrating the result over the domain $\mathbb{R}^3_+$, and using then integration by parts and \eqref{52}, we get
\begin{equation}
J(t) =2K(t).\label{62}
\end{equation}
Thus, it follows from \eqref{60}--\eqref{62} that
\begin{equation}
- 2\int_{\partial\Omega_t}\partial_t\varphi (\widehat{\mathcal{H}},\dot{\mathcal{H}})|_{x_1=0}\,{\rm d}x'\,{\rm d}s \geq K(t) - C\|\dot{\mathcal{H}}\|^2_{L_2(\Omega_t)}  \geq \frac{1}{2}\,\|\dot{\mathcal{H}}(t)\|^2_{L_2(\mathbb{R}^3_+)}
- C\|\dot{\mathcal{H}}\|^2_{L_2(\Omega_t)}.
\label{63}
\end{equation}

Consider now the term $-2(\hat{v}',\nabla'\varphi )(\widehat{\mathcal{H}},\dot{\mathcal{H}})|_{x_1=0}$ (with $\hat{v}'=(\hat{v}_2,\hat{v}_3)$) contained in \eqref{59}:
\begin{multline}
-2\int_{\partial\Omega_t}(\widehat{\mathcal{H}}_2\partial_2A+\widehat{\mathcal{H}}_3\partial_3A)
(\hat{v}',\nabla'\varphi )|_{x_1=0}\,{\rm d}x'\,{\rm d}s \\
=-2\int_{\partial\Omega_t}\Bigl\{(\hat{v}',\nabla'A)(\widehat{\mathcal{H}}',\nabla'\varphi ) \\ +
(\hat{v}_2\widehat{\mathcal{H}}_3-\hat{v}_3\widehat{\mathcal{H}}_2)(\partial_3A\partial_2\varphi -\partial_2A\partial_3\varphi )\Bigr\}\Bigr|_{x_1=0}\,{\rm d}x'\,{\rm d}s = M(t) +N(t),
\label{64}
\end{multline}
where
\[
M(t)=-2\int_{\partial\Omega_t}\dot{\mathcal{H}}_N(\hat{v}',\dot{\mathfrak{H}'})|_{x_1=0}\,{\rm d}x'\,{\rm d}s,
\]
\begin{multline*}
N(t)=-2\int_{\partial\Omega_t}\Bigl\{ \left(\partial_3 (\hat{v}_2\widehat{\mathcal{H}}_3-\hat{v}_3\widehat{\mathcal{H}}_2)+\hat{v}_2\partial_1(\widehat{\mathcal{H}}_N)\right)\dot{\mathcal{H}}_{\tau_2} \\
+ \left(\partial_2 (\hat{v}_2\widehat{\mathcal{H}}_3-\hat{v}_3\widehat{\mathcal{H}}_2)+\hat{v}_3\partial_1(\widehat{\mathcal{H}}_N)\right)\dot{\mathcal{H}}_{\tau_3}\Bigr\}\Bigr|_{x_1=0}\,\varphi\,{\rm d}x'\,{\rm d}s.
\end{multline*}
For ``frozen'' coefficients $N(t)\equiv 0$. To avoid unnecessary technical details we consider here the integral $M(t)$ for the particular case $\hat{\varphi}=0$ and leave the general case to the reader. For $\hat{\varphi}=0$
\[
M(t)=-2\int_{\partial\Omega_t}(\hat{v}_2\dot{\mathcal{H}}_1\dot{\mathcal{H}}_2 +
\hat{v}_3\dot{\mathcal{H}}_1\dot{\mathcal{H}}_3)|_{x_1=0}\,{\rm d}x'\,{\rm d}s.
\]
and the div-curl system \eqref{36} takes the form $\nabla\times\widetilde{\mathcal{H}} =0$,
${\rm div}\,\widetilde{\mathcal{H}}=0$, with $\widetilde{\mathcal{H}}=(-\dot{\mathcal{H}}_1,
\dot{\mathcal{H}}_2,\dot{\mathcal{H}}_3)$. From this system we easily deduce that
\[
\int_{\partial\Omega_t}\dot{\mathcal{H}}_1\dot{\mathcal{H}}_2|_{x_1=0}\,{\rm d}x'\,{\rm d}s=0\quad\mbox{and}\quad 
\int_{\partial\Omega_t}\dot{\mathcal{H}}_1\dot{\mathcal{H}}_3|_{x_1=0}\,{\rm d}x'\,{\rm d}s=0
\]
For the ``frozen'' coefficients case, $\hat{v}_2$ and $\hat{v}_3$ are constants and we therefore conclude that
$M(t)\equiv 0$. For variable coefficients, omitting simple calculations, from the div-curl system we derive the estimate
\begin{equation}
 M(t) \geq - C\|\dot{\mathcal{H}}\|^2_{L_2(\Omega_t)}.
\end{equation}

At last, consider the term $2\varphi(\partial_1\hat{v}_{N})(\widehat{\mathcal{H}},\dot{\mathcal{H}})|_{x_1=0}$ contained in \eqref{59}.
For variable coefficients, in view of assumptions \eqref{24} and \eqref{25}, the corresponding boundary integral together with $N(t)$ and the last integral in \eqref{60} can be written in the following compact form:
\begin{multline}
N(t)+2\int_{\partial\Omega_t}\Bigl\{(\widehat{\mathcal{H}}',\nabla'A)\,
\partial_1\hat{v}_{N} 
+(\dot{\mathfrak{H}}',\partial_t\widehat{\mathcal{H}}')\Bigr\}\Bigr|_{x_1=0}\,\varphi\,{\rm d}x'\,{\rm d}s \\
=\mathcal{N}(t)=2\int_{\partial\Omega_t}\left( \dot{\mathfrak{H}}',\bigl\{\partial_t\widehat{\mathcal{H}}-
\nabla\times (\hat{v}\times \widehat{\mathcal{H}}\,)\bigr\}'\right)\Bigr|_{x_1=0}\,\varphi\,{\rm d}x'\,{\rm d}s,
\label{66}
\end{multline}
where $\{b\}'= (b_2,b_3)$ for any three-dimensional vector $b=(b_1,b_2,b_3)$. Since we cannot control the trace $\dot{\mathcal{H}}|_{x_1=0}$ we do not know yet how to estimate the integral $\mathcal{N}(t)$ (see Section 5 for further discussions). But for ``frozen'' coefficients, taking into account \eqref{62}, we have
\begin{equation}
\mathcal{N}(t)=2(\partial_1\hat{v}_{N})\int_{\partial\Omega_t}(\widehat{\mathcal{H}}',\nabla'A)|_{x_1=0}\,
\varphi\,{\rm d}x'\,{\rm d}s= -2(\partial_1\hat{v}_{N})\int_{0}^{t}K(s)\,{\rm d}s.
\label{67}
\end{equation}

It follows from \eqref{58}--\eqref{66} that 
\begin{multline}
I(t)+\frac{1}{2}\,\|\dot{\mathcal{H}}(t)\|^2_{L_2(\mathbb{R}^3_+)}+ \frac{1}{2}\int_{\mathbb{R}^2}[\partial_1\hat{q}]\,\varphi^2\,{\rm d}x'
+\mathcal{N}(t) \\[3pt]
\leq C \Bigl\{ \| f\|^2_{L_2(\Omega_T)}
+\| g\|^2_{H^1(\partial\Omega_T)}\qquad\qquad \\[3pt]
+\int_0^t\left(I(s) +\|\dot{\mathcal{H}}(s)\|^2_{L_2(\mathbb{R}^3_+)}
+\|\varphi (s)\|^2_{L_2(\mathbb{R}^2)}\right)\,{\rm d}s\Bigr\}.\label{68}
\end{multline}
For the ``frozen'' coefficients case, by virtue of \eqref{67}, we obtain inequality \eqref{68} where we formally set $\mathcal{N}(t)=0$. Taking into account assumption \eqref{39} and applying Gronwall's lemma, from this inequality we finally deduce the basic a priori 
estimate \eqref{40}. This completes the proof of Theorem \ref{t2.2}.

\section{``Purely MHD'' case}

For the ``purely MHD'' case thanks to assumption \eqref{41} we can resolve \eqref{35}, \eqref{37}, and \eqref{44} for $\nabla_{t,x'}\varphi =(\partial_t\varphi,\nabla'\varphi)$:
\begin{equation}
\nabla_{t,x'}\varphi=\hat{a}_1\dot{H}_{N|x_1=0}+\hat{a}_2\dot{\mathcal{H}}_{N|x_1=0} +\hat{a}_3\dot{v}_{N|x_1=0}+\hat{a}_4\varphi +\hat{a}_0g_3,
\label{69}
\end{equation}
where the vector-functions $\hat{a}_{\alpha}={a}_{\alpha}(\widehat{W}_{|x_1=0})=(\hat{a}_{\alpha}^1,\hat{a}_{\alpha}^2,\hat{a}_{\alpha}^3)$ can be easily written down, in particular, $\hat{a}_{3}^2=\hat{a}_{3}^3=0$,
\[
\hat{a}_1^2=\frac{\widehat{\mathcal{H}}_3|_{x_1=0}}{(\widehat{H}_2\widehat{\mathcal{H}}_3-
\widehat{H}_3\widehat{\mathcal{H}}_2)|_{x_1=0}},\quad
\hat{a}_3^2=\frac{(\widehat{\mathcal{H}}_3\partial_1\widehat{H}_N -\widehat{H}_3\partial_1\widehat{\mathcal{H}}_N)|_{x_1=0}}{(\widehat{H}_2\widehat{\mathcal{H}}_3-
\widehat{H}_3\widehat{\mathcal{H}}_2)|_{x_1=0}},\quad \mbox{etc.}
\]
Using the terminology of paradifferential calculus, we can say that for the ``purely MHD'' case the symbol associated to the interface is {\it elliptic} (see, e.g., \cite{BS}). This fact plays the crucial role in the proof of estimate \eqref{42}.

Using the special structure of the boundary matrix $\mathcal{A}_1$ (see \eqref{50}), from system \eqref{49}, equation \eqref{43}, and estimate \eqref{47} we easily deduce the inequality
\begin{equation}
\|\partial_1V_n (t)\|^2_{L_2(\mathbb{R}^3_+)}\leq C(K)\left\{[f]^2_{2,*,T}+\nt \dot{U}(t)\nt^2_{1,*}\right\}.
\label{70}
\end{equation}
Moreover, by resolving the div-curl system \eqref{36} for the normal normal ($x_1$-) derivatives, one gets
\begin{equation}
\|\partial_1\dot{\mathcal{H}} (t)\|^2_{L_2(\mathbb{R}^3_+)}\leq C(K)\left\{\|\partial_2\dot{\mathcal{H}} (t)\|^2_{L_2(\mathbb{R}^3_+)} +\|\partial_3\dot{\mathcal{H}} (t)\|^2_{L_2(\mathbb{R}^3_+)}\right\}.
\label{71}
\end{equation}
Then, thanks to the trace theorem \eqref{70} and \eqref{71} imply
\begin{equation}
\|V_n|_{x_1=0} (t)\|^2_{L_2(\mathbb{R}^2)} \leq C \left\{[f]^2_{2,*,T}+\nt \dot{U}(t)\nt^2_{1,*}\right\},
\label{72}
\end{equation}
\begin{equation}
\|\dot{\mathcal{H}}_n|_{x_1=0} (t)\|^2_{L_2(\mathbb{R}^2)} \\ \leq C
\left\{ \|\partial_2\dot{\mathcal{H}} (t)\|^2_{L_2(\mathbb{R}^3_+)} +\|\partial_3\dot{\mathcal{H}} (t)\|^2_{L_2(\mathbb{R}^3_+)}\right\}.
\label{73}
\end{equation}

We first get an estimate for weighted derivatives. To estimate such terms we do not need boundary conditions because the weight $\sigma|_{x_1=0}=0$. By applying to system \eqref{34} the operator $\sigma\partial_1$ and using standard arguments of the energy method (see, e.g, \cite{Tr05} for more details), we obtain the inequality
\begin{equation}
\|\sigma\partial_1\dot{U}(t)\|^2_{L_2(\mathbb{R}^3_+)} \leq C\left\{ [f]^2_{1,*,T}+[\dot{U}]^2_{1,*,t} \right\}.\label{74}
\end{equation}

We can easily get the inequalities
\begin{equation}
\|\dot{U}(t)\|^2_{L_2(\mathbb{R}^3_+)}\leq C[\dot{U}]^2_{1,*,t},\quad
\|\dot{\cal H}(t)\|^2_{L_2(\mathbb{R}^3_+)}\leq C\|\dot{\cal H}\|^2_{H^1(\Omega_t)}
\label{75}
\end{equation}
and
\begin{equation}
\nt\varphi (t)\nt^2_{H^1(\mathbb{R}^2)}\leq C\Bigl\{ {[}f{]}^2_{2,*,T} +
\|g\|^2_{H^2(\partial\Omega_T)} +
 [\dot{U}]^2_{1,*,t} +\|\dot{\cal H}\|^2_{H^1(\Omega_t)}+\|\varphi\|^2_{L_2(\partial\Omega_t)}\Bigr\},
\label{76}
\end{equation}
where \eqref{75} follows from the trivial relations
\[
\frac{\rm d}{{\rm d}t}\,\|\dot{U}(t)\|^2_{L_2(\mathbb{R}^3_+)}=2\int_{\mathbb{R}^3_+}(\dot{U},\partial_t\dot{U})\,{\rm d}x,\quad
\frac{\rm d}{{\rm d}t}\,\|\dot{\cal H}(t)\|^2_{L_2(\mathbb{R}^3_+)}=2\int_{\mathbb{R}^3_+}(\dot{\cal H},\partial_t\dot{\cal H})\,{\rm d}x,
\]
and \eqref{76} is the result of the multiplication of the first boundary condition in \eqref{35} by $2\varphi$, the integration over the domain $\mathbb{R}^2$, and the usage of \eqref{69}, estimates \eqref{72} and \eqref{48}. In \eqref{76} we use the norm
\[
\nt u (t)\nt^2_{H^1(\mathbb{R}^2)}:=\int_{\mathbb{R}^2}\left( u^2 +(\partial_tu)^2  
+(\partial_2u)^2+(\partial_3u)^2\right) dx'.
\]

We now proceed to estimating the tangential derivatives $\partial_k$ and $\partial_t$ ($k=2,3$)
of the solution. This is the most important step because we shall use the boundary conditions. Differentiating system \eqref{49} with respect to $x_{\ell}$ (with $\ell=2$ or $\ell=3$) and using again standard arguments, we get the energy inequality
\begin{equation}
\int_{\mathbb{R}^3_+}({\cal A}_0\partial_{\ell}V,\partial_{\ell}V)\,{\rm d}x- 2\int_{\partial\Omega_t}\partial_{\ell}\dot{q}\,\partial_{\ell}\dot{v}_N|_{x_1=0}\,{\rm d}x'\,{\rm d}s  \leq C 
\left\{ [f]^2_{1,*,T} +[\dot{U}]^2_{1,*,t} \right\},
\label{77}
\end{equation}
By using the boundary conditions \eqref{35} we obtain
\begin{equation}
-2\int_{\partial\Omega_t}\partial_{\ell}\dot{q}\,\partial_{\ell}\dot{v}_N|_{x_1=0}\,{\rm d}x'\,{\rm d}s=
J_1(t)+J_2(t),
\label{78}
\end{equation}
where
\[
J_1(t)=2\int_{\partial\Omega_t}\Bigl\{
\left([\partial_1\hat{q}]\partial_{\ell}\varphi- \partial_{\ell}g_2\right)\partial_{\ell}\dot{v}_N
+[\partial_{\ell}\partial_1\hat{q}]\,\varphi\, \partial_{\ell}\dot{v}_N
-(\partial_{\ell}\widehat{\mathcal{H}},\dot{\mathcal{H}})\partial_{\ell}\dot{v}_N\Bigr\}\Bigr|_{x_1=0}\,{\rm d}x'\,{\rm d}s,
\]
\[ 
J_2(t)= -2\int_{\partial\Omega_t}(\widehat{\mathcal{H}},\partial_{\ell}\dot{\mathcal{H}})\partial_{\ell}\dot{v}_N|_{x_1=0}\,{\rm d}x'\,{\rm d}s.
\]

Integrating by parts and applying \eqref{69}, we get
\begin{multline*}
J_1(t)=\int_{\partial\Omega_t}\Bigl\{\hat{c}_1\dot{v}_N\partial_{\ell}\varphi +\hat{c}_2\dot{v}_N g_3 +\hat{c}_3\dot{v}_N \partial_{\ell}g_3+2\dot{v}_N\partial_{\ell}^2g_2 \\
+\hat{c}_4\dot{H}_N\partial_{\ell}\dot{v}_N+\sum\nolimits_{i=1}^{3}\hat{c}_{i+4}\dot{\cal H}_N\partial_{\ell}\dot{v}_N +\hat{c}_8\dot{v}_N\varphi
\Bigr\}\Bigr|_{x_1=0}\,{\rm d}x'\,{\rm d}s,
\end{multline*}
where 
\[
\hat{c}_{j}=\hat{c}_{j}(\widehat{W}_{|x_1=0},\partial_{\ell}\widehat{W}_{|x_1=0})\quad (j=\overline{1,7}),\quad 
\hat{c}_8=\hat{c}_8(\widehat{W}_{|x_1=0},\partial_{\ell}\widehat{W}_{|x_1=0},\partial_{\ell}^2\widehat{W}_{|x_1=0})
\]
are functions (coefficients) dependent on the basic state \eqref{21}. To treat the terms $\hat{c}_4\dot{H}_N\partial_{\ell}\dot{v}_N$ and $\hat{c}_{i+4}\dot{\cal H}_N\partial_{\ell}\dot{v}_N$ contained in the boundary integral $J_1(t)$ we use the same arguments as in \cite{Tr05,Tr}. That is, we pass to the volume integral and integrate by parts. In particular, we have:
\begin{multline}
\int_{\partial\Omega_t}\hat{c}_5\dot{H}_N\partial_{\ell}\dot{v}_N|_{x_1=0}\,{\rm d}x'\,{\rm d}s
=-\int_{\Omega_t}\partial_1\bigl(\tilde{c}_5\dot{H}_N\partial_{\ell}\dot{v}_N\bigr){\rm d}x\,{\rm d}s \\
=\int_{\Omega_t}\Bigl\{\tilde{c}_5\partial_{\ell}\dot{H}_N\partial_1\dot{v}_N +(\partial_{\ell}\tilde{c}_5)\dot{H}_N\partial_1\dot{v}_N \\-\tilde{c}_5\partial_1\dot{H}_N\partial_{\ell}\dot{v}_N-(\partial_1\tilde{c}_5)\dot{H}_N\partial_1\dot{v}_N\Bigr\}{\rm d}x\,{\rm d}s,
\label{79}
\end{multline}
where $\tilde{c}_5|_{x_1=0}=\hat{c}_5$. Taking into account \eqref{70}--\eqref{72} and \eqref{48}, we can now estimate the boundary integral $J_1(t)$ as follows:
\begin{equation}
J_1(t)\leq
C(K)\Bigl\{ {[}f{]}^2_{2,*,T} +\|g\|^2_{H^2(\partial\Omega_T)} +{[}\dot{U}{]}^2_{1,*,t}  +
\|\dot{\mathcal{H}} \|^2_{H^1(\Omega_t)}+\|\varphi\|^2_{H^1(\partial\Omega_t)} \Bigr\}.
\label{80}
\end{equation}

Regarding the integral $J_2(t)$, to treat it we can repeat arguments in \eqref{59}--\eqref{66}, where, roughly speaking, instead of $A$, $\varphi$, and $\dot{\cal H}_N$ we have $\partial_{\ell}A$, $\partial_{\ell}\varphi$, and $\partial_{\ell}\dot{\cal H}_N$ and there appear additional lower-order terms. In particular, the integral
\begin{equation}
\int_{\partial\Omega_t}\bigl({\rm coeff}',\partial_{\ell}\dot{\mathfrak{H}}'\bigr)\bigr|_{x_1=0}\,\partial_{\ell}\varphi
\,{\rm d}x'\,{\rm d}s
\label{81}
\end{equation}
is the counterpart of the integral ${\cal N}(t)$ in \eqref{66}, where ${\rm coeff}' =({\rm coeff}_2,{\rm coeff}_3)$ and ${\rm coeff}_{i}$ are coefficients dependent on the basic state \eqref{21}. But now we can treat such an integral for variable coefficients by using the ellipticity of the interface symbol. That is, we express $\partial_{\ell}\varphi$ appearing in \eqref{81} through $\dot{H}_{N|x_1=0}$, $\dot{\mathcal{H}}_{N|x_1=0}$, and $g_3$ (see \eqref{69}). After that we use the same arguments as in \eqref{79}. Omitting detailed calculations (in particular, we also exploit \eqref{73}), we finally
estimate the boundary integral $J_2(t)$:
\begin{multline}
J_2(t)\geq \frac{1}{2}\|\partial_{\ell}\dot{\cal H}(t)\|^2_{L_2(\mathbb{R}^3_+)}  -
C(K)\Bigl\{ [f]^2_{2,*,T} +\|g\|^2_{H^2(\partial\Omega_T)}\\ +[\dot{U}]^2_{1,*,t}  +
\|\dot{\mathcal{H}} \|^2_{H^1(\Omega_t)} +\|\varphi\|^2_{H^1(\partial\Omega_t)}\Bigr\}.
\label{82}
\end{multline}

It follows from \eqref{77}, \eqref{78}, \eqref{80}, and \eqref{82} that
\begin{multline}
\|\partial_{\ell}\dot{U}(t)\|^2_{L_2(\mathbb{R}^3_+)} +
\|\partial_{\ell}\dot{\cal H}(t)\|^2_{L_2(\mathbb{R}^3_+)} \leq C(K)\Bigl\{
[f]^2_{2,*,T} +\|g\|^2_{H^2(\partial\Omega_T)}\\+
[\dot{U}]^2_{1,*,t}  +
\|\dot{\mathcal{H}} \|^2_{H^1(\Omega_t)} +\|\varphi\|^2_{H^1(\partial\Omega_t)}\Bigr\}.
\label{83}
\end{multline}
Combining \eqref{83} with \eqref{71}, \eqref{74}--\eqref{76}, we get
\begin{multline}
\|\dot{U}(t)\|^2_{1,*} +
\|\dot{\cal H}(t)\|^2_{H^1(\mathbb{R}^3_+)} +\nt\varphi (t)\nt^2_{H^1(\mathbb{R}^2)} \\ \leq C(K)\Bigl\{
[f]^2_{2,*,T} +\|g\|^2_{H^2(\partial\Omega_T)} +
[\dot{U}]^2_{1,*,t}  +
\|\dot{\mathcal{H}} \|^2_{H^1(\Omega_t)} +\|\varphi\|^2_{H^1(\partial\Omega_t)}\Bigr\}.
\label{84}
\end{multline}

We still miss the $L_2$ norms of $\partial_t\dot{U}$ and $\partial_t\dot{\cal H}$ in the left-hand side of \eqref{84} which we need for closing estimate \eqref{84}. To estimate the time derivatives of $\dot{U}$ and $\dot{\cal H}$ we use the same arguments as in \eqref{77}--\eqref{83}. We just replace $\partial_{\ell}$ by $\partial_t$ (or $\partial_s$), and the only principal difference is that expressions in the form $\partial_t(\cdots )$ do not disappear 
after the integration over the domain $\Omega_t$ (whereas $\int_{\Omega_t}\partial_{\ell}(\cdots )\,{\rm d}x\,{\rm d}s=0$). For example, the right-hand side of \eqref{79} with $\partial_{\ell}$ replaced by $\partial_s$ contains the additional integral
\[
-\int_{\Omega_t}\partial_s\bigl(\tilde{c}_5\dot{H}_N\partial_{1}\dot{v}_N\bigr){\rm d}x\,{\rm d}s=-\int_{\mathbb{R}^3_+}\tilde{c}_5\dot{H}_N\partial_{1}\dot{v}_N\,{\rm d}x.
\]
Using the Young inequality, \eqref{70} , and \eqref{75}, we estimate this integral as follows:
\[
-\int_{\mathbb{R}^3_+}\tilde{c}_5\dot{H}_N\partial_{1}\dot{v}_N\,{\rm d}x\leq C(K)\left\{
[f]^2_{2,*,T}+\varepsilon\,\nt \dot{U}(t)\nt^2_{1,*} +\frac{1}{\varepsilon}\,
[\dot{U}]^2_{1,*,t}\right\},
\]
where $\varepsilon $ is a small positive constant. Dealing analogously with remaining terms appearing after the integration by parts with respect to $t$ (see also \cite{Tr05} for similar calculations), we finally obtain the inequality
\begin{multline}
\|\partial_t\dot{U}(t)\|^2_{L_2(\mathbb{R}^3_+)} +
\|\partial_t\dot{\cal H}(t)\|^2_{L_2(\mathbb{R}^3_+)}  \leq C(K)\Bigl\{
[f]^2_{2,*,T} +\|g\|^2_{H^2(\partial\Omega_T)} \\ +
[\dot{U}]^2_{1,*,t}  +
\|\dot{\mathcal{H}} \|^2_{H^1(\Omega_t)} +\|\varphi\|^2_{H^1(\partial\Omega_t)}\\+
\varepsilon\,\|\dot{U}(t)\|^2_{1,*} +\varepsilon\,
\|\dot{\cal H}(t)\|^2_{H^1(\mathbb{R}^3_+)}
\Bigr\}.
\label{85}
\end{multline}

Combining \eqref{85} with \eqref{84} and choosing the constant $\varepsilon $ to be small enough, we get 
\begin{multline}
\nt\dot{U}(t)\nt^2_{1,*} +
\nt\dot{\cal H}(t)\nt^2_{H^1(\mathbb{R}^3_+)} +\nt\varphi (t)\nt^2_{H^1(\mathbb{R}^2)} \\ \leq C(K)\Bigl\{
[f]^2_{2,*,T} +\|g\|^2_{H^2(\partial\Omega_T)} +
[\dot{U}]^2_{1,*,t}  +
\|\dot{\mathcal{H}} \|^2_{H^1(\Omega_t)} +\|\varphi\|^2_{H^1(\partial\Omega_t)}\Bigr\},
\label{86}
\end{multline}
where
\[
\nt u (t)\nt^2_{H^1(\mathbb{R}^3_+)}:=\| u (t)\|^2_{H^1(\mathbb{R}^3_+)}+
\|\partial_tu(t)\|^2_{L_2(\mathbb{R}^3_+)}.
\]
Applying Gronwall's lemma, from the last inequality we derive the desired a priori 
estimate \eqref{42}. This completes the proof of Theorem \ref{t2.3}.

\section{Concluding remarks/Open problems}

We have obtained first results towards the proof of the local-in-time existence of smooth solutions of the plasma-vacuum interface problem in ideal compressible MHD under the basic assumption that the plasma density does not go to zero continuously, but jumps. We have found two cases of well-posedness of the constant coefficient linearized problem: the ``gas dynamical'' case and the ``purely MHD'' case. From the mathematical point of view, for the ``gas dynamical'' case the interface symbol can be non-elliptic and for the ``purely MHD'' case the interface symbol is always elliptic. For the latter case the MHD counterpart \eqref{13} of the natural physical condition in gas dynamics \cite{Lind,Tr09} can be violated, i.e., the magnetic field plays a stabilizing role for well-posedness.

For the ``purely MHD'' case we have managed to derive a basic a priori estimate for the variable coefficient linearized problem. We prove this estimate in the anisotropic weighted Sobolev space $H^1_*$ because the interface is a characteristic boundary for the hyperbolic system of MHD equations, and in MHD a natural loss of derivatives in the normal direction in a priori estimates cannot be compensated as in gas dynamics \cite{Tr09}. Assuming that the original nonlinear problem has smooth enough solutions, we can easily prove the {\it uniqueness} of a solution of this problem by a standard argument and using the basic a priori estimate in $H^1_*$. 

In the basic a priori estimate \eqref{42} we have a loss of derivatives from the source terms to the solution. Clearly, we will have a loss of derivatives also in a corresponding {\it tame}
estimate whose derivation is postponed to the future. Therefore, we expect to prove the existence of solutions of the nonlinear problem by a suitable Nash-Moser-type iteration scheme.
We do not see any obstacles in this direction. We think that the derivation of the tame estimate and the realization of the Nash-Moser procedure is just a technical matter and can be done as in \cite{Tr} for current-vortex sheets. 

At the same time, there is still an open problem in getting the local-in-time existence result for the ``purely MHD'' case. The point is that we have not yet proved the existence of solutions of the linearized problem. We can naturally formulate a {\it dual problem} for it, but we still do not know how to get an a priori estimate for the dual problem. This is a very surprising fact because usually if we can obtain an a priori estimate for the original linearized problem, then in exactly the same manner we can derive it for the dual problem. After that the existence of solution of the linearized problem can be proved by the classical argument of Lax and Phillips \cite{LP}. Of course, our ``hyperbolic-elliptic'' problem is very nonstandard and this causes the mentioned difficulty. We expect to prove the existence of solutions of the linearized problem either by iterations or by considering a regularized problem. This work is postponed to the future.

Regarding the ``gas dynamical'' case, it is still unclear how to carry the basic a priori estimate obtained for the constant coefficients case over variable coefficients. The difficulty is connected with the appearance of additional lower-order terms and the fact that we cannot control the trace of the perturbation of the vacuum magnetic field in the higher norm. We think that it is very unlikely that the plasma-vacuum interface problem is not well-posed for the ``gas dynamical'' case. But, we cannot completely exclude this possibility. This question is the most important open problem both from the mathematical and the physical points of view. On the other hand, for the {\it model} free boundary problem when the vacuum magnetic field ${\cal H}\equiv 0$ we can prove a local-in-time existence theorem in the anisotropic weighted Sobolev spaces $H^m_*$, provided the initial data satisfy the condition $\partial q/\partial {\cal N}\leq -\varepsilon  <0$ (cf. \eqref{13}). Roughly speaking, we can prove such a theorem by the combination of arguments applied in \cite{Tr09} to the compressible fluid-vacuum problem and in \cite{Tr} to ideal compressible current-vortex sheets.

\vspace*{12pt}
\noindent
{\bf Acknowledgements}
\vspace*{9pt}

A part of this work was done during the two months research stay of the author at the Department of Mathematics and Statistics of the University of Konstanz sponsored by DAAD. The author gratefully thanks Heinrich Freist\"uhler
for his kind hospitality and many helpful discussions during this visit.

\end{document}